\documentclass[fleqn]{article}

\usepackage{amsmath}
\usepackage{amsfonts}
\usepackage{xcolor}
\usepackage{calrsfs}

\newcommand{\FF}{\mathbb{F}}
\newcommand{\ZZ}{\mathbb{Z}}
\newcommand{\PG}{\mathrm{PG}}

\newcommand{\cH}{{\cal H}}
\newcommand{\cS}{{\cal S}}

\newcommand{\ve}{\varepsilon}

\newcommand{\Aut}{\mathrm{Aut}}
\newcommand{\veS}{\varepsilon_{\mathrm{nat}}}
\newcommand{\veSS}{\varepsilon_{\mathrm{id}_\FF}}
\newcommand{\Pu}{\mathrm{Pure}}
\newcommand{\Ker}{\mathrm{Ker}}

\newcommand{\bx}{{\bf x}}
\newcommand{\by}{{\bf y}}

\newcommand{\be}{{\bf e}}
\newcommand{\ba}{{\bf a}}
\newcommand{\bfb}{{\bf b}}
\newcommand{\bc}{{\bf c}}
\newcommand{\bv}{{\bf v}}

\newcommand{\br}{{\bf r}}
\newcommand{\bn}{{\bf n}}
\newcommand{\bm}{{\bf m}}

\newcommand{\AnF}{A_{n,\{1,n\}}(\FF)} 
\newcommand{\idF}{\mathrm{id}_\FF} 
\newcommand{\gM}{\mathfrak{M}} 
\newcommand{\vet}{\tilde{\varepsilon}}

\newtheorem{conc}{Claim}[section]

\newtheorem{theo}[conc]{Theorem}
\newtheorem{cor}[conc]{Corollary}
\newtheorem{prop}[conc]{Proposition}
\newtheorem{lemma}[conc]{Lemma}

\newtheorem{conj}[conc]{Conjecture}
\newtheorem{note}{Remark}

\title{Embeddings and hyperplanes of the Lie geometry $A_{n,\{1,n\}}(\FF)$} 
\author{Antonio Pasini}
\date{} 

\begin{document}

\maketitle

\begin{abstract} 
In this paper we consider a family of projective embeddings of the geometry $\AnF$ of point-hyperplanes flags of $\PG(n,\FF)$. The natural embedding $\veS$ is one of them. It maps every point-hyperplane flag $(p,H)$ onto the vector-line $\langle \bx\otimes\xi\rangle$, where $\bx$ is a representative vector of $p$ and $\xi$ is a linear functional describing $H$. The other embeddings have been discovered by Thas and Van Maldeghem \cite[Part III]{TVM}  for the case $n = 2$ and later generalized to any $n$ by De Schepper, Schillewaert and Van Maldeghem \cite{DSSVM}. They are obtained as twistings of $\veS$ by non-trivial automorphisms of $\FF$. Explicitly, for $\sigma\in \Aut(\FF)\setminus\{\idF\}$, the twisting $\ve_\sigma$ of $\veS$ by $\sigma$ maps $(p,H)$ onto $\langle \bx^\sigma\otimes \xi\rangle$. We shall prove that, when $\Aut(\FF) \neq \{\idF\}$, a geometric hyperplane $\cH$ of $\AnF$ arises from $\veS$ and one of its twistings or from two distinct twistings of $\veS$ if and only if $\cH = \{(p,H)\in \AnF \mid p\in A \mbox{ or } a \in H\}$ for a possibly non-incident point-hyperplane pair $(a,A)$ of $\PG(n,\FF)$. We call these hyperplanes quasi-singular hyperplanes. With the help of this result we shall prove that if $|\Aut(\FF)| > 1$ then $\AnF$ admits no absolutely universal embedding. 
\end{abstract} 


\section{Introduction}

\subsection{Basic properties of the geometry $A_{n,\{1,n\}}(\FF)$}

Following a well established notation, we denote by $A_{n,\{1,n\}}(\FF)$ the geometry of point-hyperplane flags of the projective geometry $\PG(n,\FF)$, for $n \geq 2$ and $\FF$ a given field. 

Explicitly, $A_{n,\{1,n\}}(\FF)$ is the point-line geometry the points of which are the ordered pairs $(p,H)$ where $p$ and $H$ are a point and a hyperplane of $\PG(n,\FF)$ respectively and $p\in H$; the lines of $A_{n,\{1,n\}}(\FF)$ are the sets $\{(p,H)\mid p\in \ell\}$ for $\ell$ a line of $\PG(n,\FF)$ and $H$ a hyperplane of $\PG(n,\FF)$ containing $\ell$ and the sets $\{(p,H)\mid H \supset L\}$ for $L$ a sub-hyperplane of $\PG(n,\FF)$ (namely a subspace of $\PG(n,\FF)$ of codimension 2) and $p$ a point of $L$. Accordingly, two points $(p,H)$ and $(q,K)$ of $A_{n,\{1,n\}}(\FF)$ are collinear if and only if either $p = q$ or $H = K$. Let $p\neq q$ and $H\neq K$. If either $p\in K$ or $q\in H$, then $(p,H)$ and $(q,K)$ are at distance 2, otherwise they are at distance 3. Thus, the diameter of the collinearity graph of $A_{n,\{1,n\}}(\FF)$ is equal to 3. 




\subsubsection{Maximal singular subspaces} 

Given a point $a$ of $\PG(n,\FF)$, let ${\cal M}_a$ be the set of pairs $(a,H)$ with $H$ a hyperplane of $\PG(n,\FF)$ containing $a$. This set is a maximal singular subspace of $A_{n,\{1,n\}}(\FF)$, namely a subspace of the point-line geometry $\AnF$ all points of which are mutually collinear and maximal with respect to this property. Dually, for a hyperplane $A$ of $\PG(n,\FF)$, the set ${\cal M}_A := \{(p,A)\mid p\in A\}$ is a maximal singular subspace of $A_{n,\{1,n\}}(\FF)$. We say that ${\cal M}_a$ and ${\cal M}_A$ are {\em based} at $a$ and $A$ respectively. Every maximal singular subspace of $A_{n,\{1,n\}}(\FF)$ admits one of these two descriptions. Thus, the maximal singular subspaces of $A_{n,\{1,n\}}(\FF)$ are partioned in two families: those which are based at a point and those based at a hyperplane. Two distinct maximal singular subspaces have at most one point in common; they meet in a point only if they do not belong to the same family. Moreover. if $M$ is a maximal singular subspace of $A_{n,\{1,n\}}(\FF)$ and $(p,H)$ a point of $A_{n,\{1,n\}}(\FF)$ exterior to $\cal M$, then $(p,H)$ is collinear with at most one point of $\cal M$. 

\subsubsection{The natural embedding and its twistings}\label{veS}         

It is well known that $A_{n,\{1,n\}}(\FF)$ admits a (full) projective embedding in the projective space $\PG(M_{n+1}^0(\FF))$ of the vector space $M_{n+1}^0(\FF)$ of null-traced square matrices of order $n+1$ with entries in $\FF$, which yields the adjoint representation of the special linear group $\mathrm{SL}(n+1,\FF)$. We call it the {\em natural embedding} of $A_{n,\{1,n\}}(\FF)$ and we denote it by the symbol $\veS$. 

Explicitly, recall that $M^0_{n+1}(\FF)$ is a hyperplane of the vector space $M_{n+1}(\FF)$ of square matrices of order $n+1$ with entries in $\FF$ and the latter is the same as the tensor product $V\otimes V^*$, where $V = V(n+1,\FF)$ and $V^*$ is the dual of $V$. The pure tensors $\bx\otimes \xi$ of $V\otimes V^*$, with $\bx$ and $\xi$ non-zero vectors of $V$ and $V^*$ respectively, yield the matrices of $M_{n+1}(\FF)$ of rank $1$.  With $\bx$ and $\xi$ as above, let $[\bx]$ and $[\xi]$ be the point and the hyperplane of $\PG(n,\FF)$ represented by $\bx$ and $\xi$. Then $([\bx], [\xi])$ is a point of $A_{n,\{1,n\}}(\FF)$ if and only if $\xi(\bx) = 0$. The pure tensor $\bx\otimes \xi$, regarded as a square matrix of $M_{n+1}(\FF)$ of rank 1, is null-traced if and only if $\xi(\bx) = 0$. The natural embedding 
\[\veS :A_{n,\{1,n\}}(\FF)\rightarrow\PG(M^0_{n+1}(\FF))\]
maps the point $([\bx],[\xi])$ of $A_{n,\{1,n\}}(\FF)$ onto the point $\langle \bx\otimes \xi\rangle$ of $\PG(M^0_{n+1}(\FF))$. 

Suppose now that $\FF$ admits non-trivial automorphisms and let $\sigma$ be one of them. We can define a `twisted' version $\ve_\sigma$ of $\veS$ as follows (De Schepper, Schillewaert and Van Maldeghem \cite{DSSVM}; also Thas and Van Maldeghem \cite[Part III]{TVM}): 
\[\begin{array}{ccccc}
\ve_\sigma & : & A_{n,\{1,n\}}(\FF)& \rightarrow & \PG(M_{n+1}(\FF))\\
 & & ([\bx],[\xi]) & \rightarrow & [\bx^\sigma\otimes \xi]
\end{array}\]
where if $\bx = (x_i) = (x_i)_{i=0}^n$ then $\bx^\sigma = (x_i^\sigma)_{i=0}^n$. This mapping is indeed a projective embedding of $A_{n,\{1,n\}}(\FF)$ in $\PG(M_{n+1}(\FF))$. Note that $\dim(\ve_\sigma) = \dim(\veS) + 1$.  

\subsubsection{Quasi-singular hyperplanes of $\AnF$}\label{quasi-sing} 

Given a point $a$ and a hyperplane $A$ of $\PG(n,\FF)$, possibly $a\not\in A$, let $\cH_{a,A}$ be the set of points $(p,H)$ of $\AnF$ collinear with at  least one point of the union ${\cal M}_a\cup{\cal M}_A$ of the maximal singular subspaces ${\cal M}_a$ and ${\cal M}_A$. In other words, $(p,H)\in \cH_{a,A}$ if and only if either $p\in A$ or $a \in H$. It is proved in \cite{Pas} that $\cH_{a,A}$ is a geometric hyperplane of $\AnF$. Following \cite{Pas}, we call the hyperplanes as $\cH_{a,A}$  {\em quasi-singular} hyperplanes. When $a\in A$ then $\cH_{a,A}$ is the set of points of $\AnF$ at non-maximal distance from the point $(a,A)$. We call it a {\em singular} hyperplane.

\subsection{Embeddings and hyperplanes of point-line geometries}\label{point-line}   

In this subsection we recall some basics on projective embeddings and geometric hyperplanes of point-lines geometries. We will stick to what is necessary in order to understand the remaining parts of this introduction. A little more of information on this matter will be given in Section \ref{prel3}.   

Throughout this subsection $\Gamma$ is an arbitrary point-line geometry as defined in Shult \cite{S}, but we assume that the lines of $\Gamma$ are subsets of the set of points of $\Gamma$, henceforth denoted by $\cal P$. We also assume $\Gamma$ to be connected.

As in Shult \cite{S93}, a {\em projective embedding} of $\Gamma$ (an {\em embedding} of $\Gamma$ for short) is an injective mapping $\ve$ from the point-set $\cal P$ of $\Gamma$ to the point-set of a desarguesian projective space $\Sigma$ such that $\ve({\cal P})$ spans $\Sigma$ and $\ve(\ell) := \{\ve(p)\}_{p\in \ell}$ is a line of $\Sigma$, for every line $\ell$ of $\Gamma$. The dimension of $\Sigma$ is the {\em dimension} $\dim(\ve)$ of $\ve$.    
    
\subsubsection{Morphisms of embeddings}\label{sec emb} 

We refer to Faure and Fr\"{o}licher \cite[Definition 6.2.1]{FF} for the definition of morphisms of projective geometries. Given two embeddings $\ve_1:\Gamma\rightarrow\Sigma_1$ and $\ve_2:\Gamma\rightarrow\Sigma_2$ of $\Gamma$, a {\em morphism} from $\ve_1$ to $\ve_2$ is morphism of projective geometries $\phi:\Sigma_1\rightarrow\Sigma_2$ such that $\ve_2 = \phi\cdot\ve_1$. If $\phi$ is a morphism from $\ve_1$ to $\ve_2$ we write $\phi:\ve_1\rightarrow\ve_2$, we say that $\ve_1$ {\em covers} $\ve_2$ and we write $\ve_1\geq \ve_2$, also $\ve_2\leq \ve_1$. Note that, as $\Gamma$ is connected by assumption, if a morphism $\phi:\ve_1\rightarrow\ve_2$ exists then $\phi$ is uniquely determined by the condition $\ve_2 = \phi\cdot\ve_1$. This condition also forces $\phi:\Sigma_1\rightarrow\Sigma_2$ to be surjective. So, if $\phi$ is injective then it is an isomorphism. If this is the case then we say that $\ve_1$ and $\ve_2$ are {\em isomorphic} and we write $\ve_1\cong \ve_2$. Note that $\ve_1\cong \ve_2$ if and only if $\ve_1\geq \ve_2\geq \ve_1$.

\subsubsection{Relatively or absolutely universal embeddings} 

Following Shult \cite{S93}, we say that an embedding $\ve$ of a point-line geometry $\Gamma$ is {\em relatively universal} if every morphism from an embedding of $\Gamma$ to $\ve$ is an isomorphism. An embedding $\ve$ is said to be {\em absolutely universal} if it covers all embeddings of $\Gamma$. Clearly, the absolutely universal embedding, if it exits, is unique up to isomorphisms and it is relatively universal. The converse is false in general. Indeed, as proved by Ronan \cite{Ron}, every embedding $\ve$ of a point-line geometry $\Gamma$ is covered by a relatively universal embedding, uniquely determined by $\ve$ up to isomorphisms. So, $\Gamma$ admits the absolutely universal embedding if and only if it admits a unique relatively universal embedding (unique
up to isomorphisms, of course).  

\begin{note}
\em
The previous definition of the absolutely universal embedding is more restrictive than in \cite{S93}, where an embedding of $\Gamma$ is said to be absolutely universal if it covers all embeddings of $\Gamma$ defined over the same division ring as $\ve$ itself. If $\Gamma$ admits embeddings defined over different division rings then it cannot admit the absolutely universal embedding as we have defined it but some of its embeddings could be absolutely universal in the sense of Shult \cite{S93}. 
\end{note}        

\subsubsection{Hyperplanes and embeddings}\label{sec hyp}

A {\em subspace} of a point-line geometry $\Gamma$ is a subset $X$ of the point-set $\cal P$ of $\Gamma$ such that, for every line $\ell$ of $\Gamma$, if $|\ell\cap X| > 1$ then $\ell\subseteq X$. A proper subspace of $\Gamma$ is said to be a {\em geometric hyperplane} of $\Gamma$ (a {\em hyperplane} of $\Gamma$ for short) if every line of $\Gamma$ meets it non-trivially. Equivalently, a hyperplane of $\Gamma$ is a proper subset $\cH$ of $\cal P$ such that, for every line $\ell$ of $\Gamma$, either $\ell \subseteq \cH$ or $|\ell\cap \cH| = 1$. 

Assuming that $\Gamma$ is embeddable, let $\ve:\Gamma\rightarrow\Sigma$ be an embedding of $\Gamma$ in a projective space $\Sigma$. We say that a hyperplane $\cH$ of $\Gamma$ {\em arises from} $\ve$ if $\ve(\cH)$ spans a projective hyperplane $H$ of $\Sigma$ and $\cH = \ve^{-1}(H)$. Conversely, if $H$ is a projective hyperplane of $\Sigma$ then $\ve^{-1}(H)$ is a hyperplane of $\Gamma$ and $\ve(\ve^{-1}(H)) = H\cap \ve({\cal P})$, but in general $\ve^{-1}(H)$ does not arise from $\ve$. Indeed $\ve^{-1}(H)$ arises from $\ve$ if and only if $\langle H\cap\ve({\cal P})\rangle = H$. The following is proved in \cite[Proposition 1.1]{Pas}. 

\begin{prop}\label{Gen 1} 
Given an embedding $\ve:\Gamma\rightarrow\Sigma$ of $\Gamma$ and a hyperplane $\cH$ of $\Gamma$, if $\cH$ is a maximal subspace of $\Gamma$ then either $\cH$ arises from $\ve$ or $\ve(\cH)$ spans $\Sigma$. In particular, given a projective hyperplane $H$ of $\Sigma$, suppose that $\ve^{-1}(H)$ is a maximal subspace of $\Gamma$. Then $H\cap\ve({\cal P})$ spans $H$. 
\end{prop}

\begin{note}
\em
A subspace $S$ of $\Gamma$ is said to {\em arise from} an embedding $\ve:\Gamma\rightarrow\Sigma$ if $S = \ve^{-1}(\langle \ve(S)\rangle)$. A hyperplane $\cH$ of $\Gamma$ arises from $\ve$ as a hyperplane precisely when it arises from $\ve$ as a subspace of $\Gamma$ and $\ve(\cH)$ spans a hyperplane of $\Sigma$. 
\end{note} 

\subsubsection{Generating rank and relatively universal embeddings}\label{gen rk}

The intersection of an arbitrary family of subspaces of $\Gamma$ is still a subspace of $\Gamma$. For a subset $X$ of $\cal P$, let $\langle X\rangle$ be the smallest subspace of $\Gamma$ containing $X$, namely $\langle X\rangle$ is the intersection of all subspaces of $\Gamma$ which contain $X$. If $\langle X\rangle = \Gamma$ then we say that $X$ {\em generates} $\Gamma$. The {\em generating rank} $\mathrm{grk}(\Gamma)$ of $\Gamma$ is the minumum  cadinality of a generating set of $\Gamma$. Clearly, if $\ve$ is a projective embedding of $\Gamma$ then $\dim(\ve) + 1 \leq \mathrm{grk}(\Gamma)$. Consequently, 

\begin{prop}\label{standard arg} 
If $\dim(\ve) + 1 = \mathrm{grk}(\Gamma) < \infty$ then $\ve$ is relatively universal. 
\end{prop}  

\subsection{Back to $\AnF$ and its embeddings}

As we shall prove in this paper (Corollary \ref{Main5}), when $|\mathrm{Aut}(\FF)| > 1$ the geometry $\AnF$ admits no absolutely universal embedding. 

\begin{conj}\label{Conj Main}
If $|\mathrm{Aut}(\FF)| = 1$ then $\veS$ is absolutely universal.  
\end{conj}

The following result by Blok and Pasini \cite{BP1} is a clue in favour of Conjecture \ref{Conj Main}: if $\FF$ is a prime field and $n > 2$ then $\veS$ is absolutely universal.  

\begin{conj}\label{Main Conj}
Both $\veS$ and all of its twistings are relatively universal.
\end{conj}

Partial results are known which support Conjecture \ref{Main Conj}. For instance, it follows from V\"{o}lklein \cite{V} that if $\FF$ is perfect of positive characteristic or an algebraic number field then $\veS$ is relatively universal. Moreover, let $\FF$ be finite. Then the generating rank of $\AnF$ is at most $(n+1)^2$ (Blok and Pasini \cite{BP2}). The twistings of $\veS$ have dimension equal to $(n+1)^2-1$. Hence all of them are relatively universal by Proposition \ref{standard arg}. As we shall notice in the last paragraph of Section \ref{sec last}, when $\FF$ is finite and $n = 2$ the same conclusion is implicit in the classification of polarized embeddings of $A_{2,\{1,2\}}(\FF)$ with $\FF$ finite by Thas and Van Maldeghem \cite{TVM}.  

Turning to the hyperplanes of $\AnF$, it is proved in \cite[Theorem 1.5]{Pas} that all hyperplanes of $\AnF$ are maximal subspaces. Therefore, in view of Proposition \ref{Gen 1}, 

\begin{prop}\label{Main0} 
For every projective embedding $\ve:\AnF\rightarrow\Sigma$, the hyperplanes of $\AnF$ which arise from $\ve$ are precisely the $\ve$-preimages of the hyperplanes of $\Sigma$.
\end{prop} 
 
\subsection{Main results}\label{Main sec} 

In view of the next theorem it is convenient to slightly modify our notation. The natural embedding $\veS$ of $\AnF$ can be regarded as a borderline case of twisting, namely the twisting of $\veS$ by the identity automorphism $\mathrm{id}_\FF$ of $\FF$. Accordingly, henceforth we will also denote $\veS$ by the symbol $\veSS$. With this notation, $\veS = \veSS$ is a member of the set $\{\ve_\sigma\}_{\sigma\in \Aut(\FF)}$. 

\begin{theo}\label{Main1}
Suppose that $|\mathrm{Aut}(\FF)| > 1$ and let $\cH$ be a geometric hyperplane of $\AnF$. 
\begin{itemize}
\item[$(1)$] If $\cH$ is quasi-singular then $\cH$ arises from $\ve_\sigma$ for every $\sigma\in\Aut(\FF)$.
\item[$(2)$] If $\cH$ is not quasi-singular then $\cH$ arises from $\ve_\sigma$ for at most one $\sigma\in\Aut(\FF)$. 
\end{itemize}
\end{theo}

Claim (1) of Theorem \ref{Main1} will be proved at the end of Section \ref{tensor}. Claim (2) will be proved in Section \ref{Proof}.  

As we shall see in Section \ref{tensor}, for every $\sigma\in \Aut(\FF)$ many hyperplanes of $\AnF$ exist which arise from $\ve_\sigma$ but are not quasi-singular. Therefore Theorem \ref{Main1} implies the following. 

\begin{cor}\label{Main2}
Let $|\mathrm{Aut}(\FF)|  > 1$. Then $\ve_\sigma\not\cong\ve_\rho$ for any choice of distinct automorphism $\sigma$ and $\rho$ of $\FF$. 
\end{cor}

The following is the second main result of this paper. We shall prove it in Section \ref{Absolute} with the help of Theorem \ref{Main1}.

\begin{theo}\label{Main4}
For any two distinct autmorhisms $\sigma$ and $\rho$ of $\FF$, no projective embedding of $\AnF$ covers both $\ve_\sigma$ and $\ve_\rho$. 
\end{theo}    

The next corollary immediately follows from Theorem \ref{Main4}. 

\begin{cor}\label{Main5} 
If $|\mathrm{Aut}(\FF)| > 1$ then $\AnF$ admits no absolutely universal embedding. 
\end{cor} 

\noindent
{\bf Organization of the paper.} In Section \ref{prel} we fix some notation for vectors, tensors, matrices and projective points, we recall some basics on semi-polynomials and we add some information on embeddings and hyperplanes of point-line geometries, as a completion of Section \ref{point-line}. 

In Section \ref{tensor} we describe the hyperplanes of $\AnF$ which arise from $\ve_\sigma$ for $\sigma\in \Aut(\FF)$. We will do it with the help of a symmetric bilinear form, called the {\em saturation form} in \cite{Pas}. We have exploited that form in \cite[Section 1.3.1]{Pas} in order to describe the hyperplanes of $\AnF$ which arise from $\veS$. Section \ref{tensor} is indeed a survey of Sections 1.3.1 and 2.2 of \cite{Pas}, with the addition of a straightforward generalization to $\ve_\sigma$ for $\sigma\neq \mathrm{id}_\FF$. 

A proof of claim (1) of Theorem \ref{Main1} will also be given at the end of Section \ref{tensor}. Claim (2) of Theorem \ref{Main1} will be proved in Section \ref{Proof}. Section \ref{Absolute} is devoted to the proof of Theorem \ref{Main4}. 

In the last section of this paper we shall consider quotients (namely morphic images) of $\veS$ and its twistings, focusing on the polarized ones, an embedding of $\AnF$ being called {\em polarized} when all singular hyperplanes of $\AnF$ arise from it. In particular, we shall generalize a part of the main result of the triad of papers \cite{TVM} by J. A. Thas and H. Van Maldeghem. 

\section{Preliminaries}\label{prel}

\subsection{Notation to be used throughout this paper} 

We denote the vectors of $V = V(n+1,\FF)$ by low case boldface roman letters and the vectors of its dual $V^*$ by low case greek letters, but we denote both the null-vector of $V$ and the null-vector of $V^*$ by the symbol $0$. Most of the matrices to be considered in this paper are square matrices of order $n+1$. We denote them by capital roman letters. In particular, $O$ and $I$ are the null square matrix and the identity matrix of order $n+1$, respectively. Scalars will be denoted by roman or greek low case letters.  

The symbols $M_{n+1}(\FF)$ and $M_{n+1}^0(\FF)$ denote respectively the vector space of all square matrices of order $n+1$ with entries in $\FF$ and the hyperplane of $M_{n+1}(\FF)$ formed by the null-traced matrices, as in Section \ref{veS}. Recall that $M_{n+1}(\FF)$ is basically the same as $V\otimes V^*$. Accordingly, we will freely switch from the matrix notation to the tensor notation and conversely, whenever these changes of notation will be convenient.  

Given a basis $E = (\be_i)_{i=0}^n$ of $V$ let $E^* = (\eta_i)_{i=0}^n$ the basis of $V^*$ dual to it. So, $\eta_i(\be_j) = \delta_{i,j}$ (Kronecker symbol) for any choice of $i, j = 0, 1,..., n$. The pure tensors $E_{i,j} = \be_i\otimes \eta_j$ form a basis of $V\otimes V^*$, henceforth denoted by $E\otimes E^*$. When $E$ is the natural basis of $V = V(n+1,\FF)$ then $E\otimes E^*$ is the natural basis of $V\otimes V^*$ and $E_{i,j}$, regarded as a matrix, is the matrix where all entries are null but the $(i,j)$-entry, which is $1$. So, $E\otimes E^*$ is just the usual natural basis of $M_{n+1}(\FF)$.

Henceforth we always assume that $E = (\be_i)_{i=0}^n$ is the natural basis of $V$ and $(\eta_i)_{i=0}^n$ is the dual $E^*$ of $E$. Accordingly, a vector $\bx = \sum_{i=0}\be_ix_i$ of $V$ is the same as the $(n+1)$-tuple $(x_i)_{i=0}^n$ of its coordinates with respect to $E$, this $(n+1)$-tuple being regarded as a column, namely an $n\times 1$ matrix. Similarly, every vector $\xi = \sum_{i=0}\xi_i\eta_i \in V^*$ is the same as the $1\times n$ matrix $(\xi_0, \xi_1,..., \xi_n)$. Thus, the scalar $\xi(\bx)$ is the same as the row-times-column product $\xi\bx$ and, for a matrix $M\in M_{n+1}(\FF)$, the product $\xi M\bx$ is the product of the row $\xi$ times $M$ times the column $\bx$. The tensor $\bx\otimes\xi$ is the same as the column-times-row product $\bx\xi$.  

We denote by $\Pu(V\otimes V^*)$ ($= \Pu(M_{n+1}(\FF))$) the set of pure tensors of $V\otimes V^* = M_{n+1}(\FF)$, with the convention that the null vector $O$ of $V\otimes V^*$ is not included in that set. So, $\Pu(M_{n+1}(\FF))$ is the set of matrices of $M_{n+1}(\FF)$ of rank 1. The hyperplane $M^0_{n+1}(\FF)$ of $M_{n+1}(\FF)$ will also be denoted by $(V\otimes V^*)_0$ and we set $\Pu((V\otimes V^*)_0) := \Pu(V\otimes V^*)\cap(V\otimes V^*)_0$. So, $\Pu(M_{n+1}^0(\FF))$ ($= \Pu((V\otimes V^*)_0)$) is the set of null-traced matrices of $M_{n+1}(\FF)$ of rank $1$, namely the set of pure tensors $\bx\otimes \xi$ with $\xi(\bx) = 0$.

Given two non-zero vectors $\bx, \by \in V$, if $\bx$ and $\by$ are proportional we write $\bx \equiv \by$. We use the same notation for vectors of $V^*$ and matrices of $M_{n+1}(\FF)$. For instance, for a matrix $M\in M_{n+1}(\FF)\setminus\{O\}$, when writing $M \equiv I$ we mean that $M$ is a scalar matrix. 

When we need to distinguish betwen a non-zero vector $\bx$ of $V$ and the point of $\PG(V) = \PG(n,\FF)$ represented by it, we denote the latter by $[\bx]$. We extend this convention to subsets of $V$: if $X\subseteq V\setminus\{0\}$ then $[X] := \{[\bx]\mid \bx\in X\}$. The same conventions will be adopted for vectors and subsets of $V^*$ and $V\otimes V^*$. In particular, if $\xi\in V^*\setminus\{0\}$ then $[\xi]$ is the point of $\PG(V^*)$ which corresponds to the hyperplane $[\Ker(\xi)]$ of $\PG(V)$. In the sequel we shall freely take $[\xi]$ as a name for $[\Ker(\xi)]$. 

\subsection{Semi-polynomials}

Let $\mathfrak{M}(\FF)$ be the module over the ring $\ZZ$ of integers with $\Aut(\FF)$ as a basis, namely the elements of $\mathfrak{M}(\FF)$ are finite formal combinations $\sum_{\sigma\in\Aut(\FF)}k_\sigma\sigma$ with $k_\sigma\in \ZZ$ for every $\sigma\in\Aut(\FF)$ and $k_\sigma = 0$ for all but at most a finite number of choices of $\sigma\in \Aut(\FF)$. Put
\[\mathfrak{M}^+(\FF) ~ := ~ \{\sum_{\sigma\in\Aut(\FF)}k_\sigma\sigma \mid k_\sigma \geq 0 \mbox{ for every } \sigma\in \Aut(\FF)\}.\]
The sum $w(\gamma) := \sum_{\sigma\in\Aut(\FF)}k_\sigma$ will be called the {\em weight} of an element $\gamma = \sum_{\sigma\in\Aut(\FF)}k_\sigma\sigma$ of $\mathfrak{M}^+(\FF)$. Given $\gamma = \sum_{i=1}^rk_i\sigma_i \in \mathfrak{M}^+(\FF)$ and $t\in \FF$, the exponentiation $t^\gamma$ stands for $\prod_{i=1}^r(t^{k_i})^{\sigma_i}$. 

Note that, regarded $\mathfrak{M}(\FF)$ as an additive monoid, $\mathfrak{M}^+(\FF)$ is a submonoid of $\mathfrak{M}(\FF)$. However $\mathfrak{M}^+(\FF)$ is not yet the structure we need. Indeed it can happen that $\mathfrak{M}^+(\FF)$ contains distinct elements $\gamma_1$ and $\gamma_2$ such that $t^{\gamma_1} = t^{\gamma_2}$ for every $t\in \FF$. (This is indeed always the case when $\FF$ is finite.) In view of what we are going to do in the sequel, this possibility would cause some troubles. We need to counteract its effects.   

So, let $\sim$ be the relation defined on $\gM^+(\FF)$ by the following clause: $\gamma_1 \sim \gamma_2$ if and only if $t^{\gamma_1} = t^{\gamma_2}$ for every $t \in \FF$. Then $\sim$ is a congruence relation of $\gM^+(\FF)$, namely it is an equivalence relation and it is preserved when taking sums. The quotient $\gM^{[+]}(\FF) := \gM^+(\FF)/\sim$ is the structure we need. The elements of $\gM^{[+]}(\FF)$ are the equivalence classes of $\sim$. For $\gamma\in \gM^+(\FF)$, we denote by $[\gamma]$ the class of $\sim$ which contains $\gamma$. The sum of $\gM^{[+]}(\FF)$ is defined as follows: $[\gamma_1]+[\gamma_2] := [\gamma_1+\gamma_2]$, for any two elements $\gamma_1, \gamma_2$ of $\gM^+(\FF)$. The {\em weight} $w(X)$ of a class $X$ of $\sim$ is the minimum among the weights of its elements. For instance, if $X$ contains $k\cdot\idF$ for $0 \leq k < |\FF|$ then $k = w(X)$. Note that if $\FF$ is finite then $w(X) < |\FF|$, for every $X\in \gM^{[+]}(\FF)$.

We are now ready to define semi-monomials and semi-polynomials. A (non-null) {\em semi-monomial} ({\em over $\FF$}) in the unknowns $t_1,..., t_m$ is a formal product 
\[M ~ = ~ a\cdot t_1^{X_1}...t_m^{X_m}\]
where $X_1,..., X_m\in \mathfrak{M}^{[+]}(\FF)$ and $a \in \FF\setminus\{0\}$. Its {\em degree} $\mathrm{deg}(M)$ is the sum $\sum_{i=1}^mw(X_i)$ of the exponents $X_1,..., X_m$ and its {\em type} is the $m$-tuple $(X_1,..., X_m)$. We allow $X = 0$ (the null element of $\gM^{[+]}(\FF)$) and we put $t^0 := 1$, as usual. We also admit the null monomial as a semi-monomial in whatever set of unknown we like, giving it $-\infty$ as its degree and $\emptyset$ as its type.     

A (non-null) {\em semi-polynomial} ({\em over $\FF$}) in the unknowns $t_1,..., t_m$ is a formal sum $P = M_1+...+M_s$ of non-null semi-monomials (over $\FF$) in those unknowns, no two of which have the same type. The {\em degree} $\mathrm{deg}(P)$ of $P$ is defined as the maximum $\mathrm{max}_{i=1}^s\mathrm{deg}(M_i)$. We allow $P$ to be the null semi-monomial, henceforth also called the {\em null semi-polynomial}. 

The following statement will be freely used troughout Section \ref{sec riciclo}. We believe this statement is well known. However, as we couldn't find any mention of it in the literature, we shall give of its proof here, so that the reader can save the trouble of looking for it in the literature.

\begin{theo}[Identity Principle]\label{Princ}
Let $P = P(t_1,..., t_m)$ be a non-null semi-polynomial over $\FF$ in the unknowns $t_1,..., t_m$. Then $P(c_1,..., c_m) \neq 0$ for at least one $m$-tuple $(c_1,..., c_m) \in \FF^m$. 
\end{theo}
{\bf Proof.} The proof exploits a double induction, on the number $m$ of unknowns of $P$ and the number $s$ of semi-monomials of $P$. Let $m = 1$. If $s = 1$ there is nothing to prove. Suppose $s = 2$. So $P = P(t) = at^X + bt^Y$ with $X, Y$ distinct elements of $\gM^{[+]}(\FF)$ and $a,b \in \FF\setminus\{0\}$. By way of contradiction, suppose that $P(c) = 0$ for every $c\in \FF$. Then with $c = 1$ we get $ a + b = 0$. It follows that $c^X = c^Y$ for every $c\in \FF$. Hence $X = Y$ by the definition of $\gM^{[+]}(\FF)$, contradicting the assumption that no two semi-monomials of the same semi-polynomial have the same type. This fixes the case $m = 1$ and $s = 2$. 

Still with $t = 1$, suppose now that $s > 2$. Now $P(t) = a_1t^{X_1}+... + a_st^{X_s}$ where $a_i \neq 0$ for every $i = 1,..., s$ and $X_i \neq X_j$ for any choice of $i, j = 1, 2,..., n$ with $i \neq j$. With no loss, we can assume that $a_s = 1$. For a contradiction, suppose that $P(c) = 0$ for every $c\in \FF$. Then $c^{X_s} = -\sum_{i=1}^{s-1}a_ic^{X_i}$ for every $c\in \FF$. Therefore, for every choice of $c, d\in\FF$ we have 
\[(\sum_{i=1}^{s-1}a_ic^{X_i})(\sum_{i=1}^{s-1}a_id^{X_i}) ~ = ~ c^{X_s}d^{X_s} ~ = (cd)^{X_s} = -\sum_{i=1}^{s-1}a_i(cd)^{X_i}.\]
Put $f(d) := \sum_{i=1}^{s-1}a_id^{X_i}$. Then, for every choice of $d$, we have $a_if(d) = -a_id^{X_i}$ by the inductive hypothesis on the polynomial
$\sum_{i=1}^sa_i(f(d)+d^{X_i})t^{X_i}$. Equivalently, since $a_i\neq 0$ by assumption, $f(d) = -d^{X_i}$. This holds for every $d\in \FF$ and every $i = 1,..., s-1$. Accordingly, $d^{X_i} = d^{X_j}$ for every $d\in \FF$ and every $i, j = 1,..., s-1$ with $i \neq j$. This forces $X_i = X_j$ by the definition of $\gM^{[+]}(\FF)$, contradicting the assumption that $X_i \neq X_j$ if $i \neq j$. By induction, the statement of the theorem holds true when $m = 1$.

Let now $m > 1$. Then $P(t_1,..., t_m) = \sum_{j=1}^rP_j(t_1,..., t_{m-1})t_m^{X_j}$ for suitable polynomials $P_1,..., P_r$ in the unknowns $t_1,..., t_{m-1}$ and suitable pairwise distinct elements $X_1,..., X_r \in \gM^{[+]}(\FF)$. 

For a contradiction, suppose that $P(c_1,..., c_{m-1}, c_m) = 0$ for every choice of $(c_1,..., c_m)\in \FF^m$. For every $(m-1)$-tuple $\bc = (c_i)_{i=1}^{m-1}\in \FF^{m-1}$, put $P_\bc(t) = P(c_1,..., c_{m-1}, t)$. Then $P_\bc(c) = 0$ for every $c\in \FF$. Hence $P_j(c_1,..., c_{m-1}) = 0$ for every $j = 1,..., r$, by the previous paragraph. However $\bc$ is an arbitrary $(m-1)$-tuple in $\FF^{m-1}$. Therefore, by the inductive hypothesis, $P_j$ is null for every $j = 1,..., r$. Hence $P$ is null as well, contradicting the assumption that $P$ is not the null polynomial.  \hfill $\Box$ 

\subsection{Addendum to Section \ref{point-line}}\label{prel3} 

As in Section \ref{point-line}, throughout this subsection $\Gamma$ is a connected point-line geometry and $\cal P$ is its point-set.   

\subsubsection{Morphisms and quotients of embeddings}\label{prel3 sub} 

Given two projective embeddings $\ve_1:\Gamma\rightarrow\Sigma_1$ and $\ve_2:\Gamma\rightarrow\Sigma_2$, suppose that $\ve_1\geq \ve_2$ and let $\phi:\Sigma_1\rightarrow\Sigma_2$ be the morphism from $\ve_1$ to $\ve_2$. Let $K := \mathrm{Ker}(\phi)$ be the kernel of $\phi$ (defined as in \cite[Definition 6.1.1]{FF}) and $\pi_K$ the projection of $\Sigma_1$ onto the star $\Sigma_1/K$ of $K$ in $\Sigma_1$. Then $K$ satisfies the following property 
\begin{equation}\label{eq quot}
K \cap \langle \ve_1(p), \ve_1(q)\rangle ~ = ~ \emptyset, \hspace{5 mm} \forall p, q \in {\cal P}
\end{equation}
and the mapping $\ve_1/K := \pi_K\cdot\ve_1$ is an embedding of $\Gamma$ isomorphic to $\ve_2$. We call $\ve_1/K$ the {\em quotient} of $\ve_1$ over $K$. Following a well established custom, we also say that $\ve_2$ is a {\em quotient} of $\ve_1$, as if $\ve_2$ was the same as $\ve_1/K$. Admittedly, this is an abuse, but it is harmless. 

Note that (\ref{eq quot}) characterizes the subspaces $K$ of $\Sigma$ such that the mapping $\ve_1/K = \pi_K\cdot\ve_1$ is an embedding of $\Gamma$. We say that a subspace $K$ of $\Sigma_1$ {\em defines a quotient} of $\ve_1$ if it satisfies (\ref{eq quot}).  

\subsubsection{More information on hyperplanes and embeddings}  

Let $\ve_1$ and $\ve_2$ be two projective embeddings of $\Gamma$. Clearly, if $\ve_1\cong \ve_2$ then a hyperplane of $\Gamma$ arises from $\ve_1$ if and only if it arises from $\ve_2$. More generally the following holds.   

\begin{prop}\label{Gen 2}
Given a projective embedding $\ve:\Gamma\rightarrow\Sigma$ let $K$ be a subpace of $\Sigma$ defining a quotient of $\ve$. Then all hyperplanes of $\Gamma$ which arise from $\ve/K$ and are maximal as subspaces of $\Gamma$ also arise from $\ve$. 
\end{prop}
{\bf Proof.} Let $\cH$ be a hyperplane of $\Gamma$ such that $\cH$ is a maximal subspace of $\Gamma$ and its $\ve/K$-image $(\ve/K)(\cH)$ spans a projective hyperplane $H$ of $\Sigma/K$. As $\cH$ is a maximal subspace of $\Gamma$, the $\ve/K$-preimage of $H$ is equal to $\cH$. The preimage $H' := \pi^{-1}_K(H) = \langle \ve(\cH)\cup K\rangle$ of $H$ by the projection $\pi_K$ of $\Sigma$ onto $\Sigma/K$ is a hyperplane of $\Sigma$ and $\cH = \ve^{-1}(H')$. It remains to prove that $H' = \langle \ve(\cH)\rangle$. For a contradiction, suppose that $\langle\ve(\cH)\rangle \subset H'$. Hence $\langle \ve(\cH\cup\{p\})\rangle \subset \Sigma$ for every point $p$ of $\Gamma$. This holds even if $p\not\in\cH$. However, if $p\not\in\cH$ then $\cH\cup\{p\}$ generates $\Gamma$, by the maximality of $\cH$. So, $\ve({\cal P})$ does not span $\Sigma$. This contradicts the definition of embedding. Therefore $H' = \langle \ve(\cH)\rangle$.  \hfill $\Box$   

\begin{prop}\label{quot0}
Given a projective embedding $\ve:\Gamma\rightarrow \Sigma$ of $\Gamma$ and a subspace of $\Sigma$, suppose that $K$ defines a quotient of $\ve$. A hyperplane $\cH$ of $\Gamma$ arises from $\ve/K$ if and only if it arises from $\ve$ and the span of $\ve(\cH)$ in $\Sigma$ contains $K$. 
\end{prop}
{\bf Proof.} The `only if' part is implicit in the proof of Proposition \ref{Gen 2}. The 'if' part is obvious. \hfill $\Box$  
 
\section{The saturation form}\label{tensor}  

Let $f : M_{n+1}(\FF)\times M_{n+1}(\FF) \rightarrow \FF$ be the symmetric bilinear form defined as follows: 
\begin{equation}\label{f1}
f(X, Y) ~ := ~ \mathrm{trace}(XY) ~ \mbox{ for any two matrices } X,Y \in M_{n+1}(\FF), 
\end{equation}
where $XY$ is the usual row-times-column product. So, with $X = (x_{i,j})_{i,j=0}^n$ and $Y = (y_{i,j})_{i,j=0}^n$ we have  
\begin{equation}\label{f2} 
f((x_{i,j})_{i,j=0}^n, (y_{i,j})_{i,j=0}^n) ~ = ~ \sum_{i,j}x_{i,j}y_{j,i}.
\end{equation}
In particular, when $X$ and $Y$ have rank 1, namely they are pure tensors, say $X = \bx\otimes\xi$ and $Y = \by\otimes \upsilon$, then
\begin{equation}\label{f3} 
f(\bx\otimes\xi, \by\otimes\upsilon) ~ = ~ \upsilon(\bx)\xi(\by).
\end{equation}
Following \cite{Pas}, we call $f$ the {\em saturation form} of $M_{n+1}(\FF)$. 

As noticed in \cite[Section 2.2]{Pas}, the form $f$ is non-degenerate. We denote by $\cS_f$ the set of $f$-isotropic vectors of $M_{n+1}(\FF)$, namely the set of matrices $X\in M_{n+1}(\FF)$ such that $f(X,X) = 0$.  

If $\mathrm{rank}(X) = 1$, namely $X = \bx\otimes \xi \in \Pu(V\otimes V^*)$, then $f(X,X) = (\xi(\bx))^2$ by formula (\ref{f3}) of Section \ref{tensor}.  Hence $f(X,X)= 0$ if and only if $\xi(\bx) = 0$, namely $X\in M_{n+1}^0(\FF)$. Accordingly,  
\[\cS_f\cap \Pu(V\otimes V^*) ~ = ~ \Pu((V\otimes V^*)_0).\] 
$\Pu((V\otimes V^*)_0)$ spans $(V\otimes V^*)_0$. Note also that $[\Pu((V\otimes V^*)_0)] = \veS({\cal P})$, where $\cal P$ stands for the set of points of $A_{n,\{1,n\}}(\FF)$. On the other hand, let $\sigma \in \Aut(\FF)\setminus\{\mathrm{id}_\FF\}$. Then $\ve_\sigma({\cal P}) = [P_\sigma]$ where \[P_\sigma ~:= ~\{\bx^\sigma\otimes \xi \mid\bx\otimes\xi \in \Pu((V\otimes V^*)_0)\}.\]
The set $P_\sigma$ contains $\be_i\otimes\eta_j$ and $(\be_i+\be_j)\otimes(\eta_i-\eta_j)$ for every choice of $i \neq j$. These vectors span $(V\otimes V^*)_0$. However $P_\sigma$ also contains pure tensors exterior to $(V\otimes V^*)_0$. For instance, $(\be_0t^\sigma + \be_1)\otimes(\eta_0-t\eta_1)$ with $t\in \FF$ such that $t^\sigma \neq t$ is one of them. Therefore $P_\sigma$ spans $V\otimes V^*$, namely $\ve_\sigma({\cal P})$ spans $\PG(M_{n+1}(\FF))$. 

Let $\perp_f$ be the orthogonality relation associated to $f$. As $f$ is non-degenerate, the hyperplanes of $M_{n+1}(\FF)$ are the perps $M^{\perp_f}$ for $M\in M_{n+1}(\FF)\setminus\{O\}$ and, for two matrices $M, N\in M_{n+1}(\FF)\setminus\{O\}$, we have $M^{\perp_f} = N^{\perp_f}$ if and only if $M \equiv N$. It is clear from formula (\ref{f1}) that  $I^{\perp_f} = M_{n+1}^0(\FF)$. Therefore, for $M\in M_{n+1}(\FF)\setminus\{O\}$, we have $M^{\perp_f} = M_{n+1}^0(\FF)$ if and only if $M\equiv I$.

Every hyperplane of $M_{n+1}^0(\FF)$ is the intersection of $M^0_{n+1}(\FF)$ with a hyperplane $M^{\perp_f}$ of $M_{n+1}(\FF)$ for a suitable matrix $M \not\in \langle I\rangle$. Of course the equality $M^{\perp_f}\cap M_{n+1}^0(\FF) = N^{\perp_f}\cap M_{n+1}^0(\FF)$ does not imply that $M\equiv N$. The next statement follows from well known properties of the polarities associated to non-degenerate reflexive bilinear forms.  

\begin{prop}\label{prop0}
For $M, N\in M_{n+1}(\FF)\setminus\{O\}$, we have $M^{\perp_f}\cap M_{n+1}^0(\FF) = N^{\perp_f}\cap M_{n+1}^0(\FF)$ if and only if 
$\langle M, I\rangle = \langle N, I\rangle$.  
\end{prop}

When $X \in \Pu(V\otimes V^*)$, namely $X\in M_{n+1}(\FF)$ has rank equal to 1, the orthogonality condition $X\perp_f M$ admits an easy formulation. Indeed formula (\ref{f2}) yields the following.

\begin{prop}\label{prop1}
Let $\bx\in V\setminus\{0\}$, $\xi\in V^*\setminus\{0\}$ and $M\in M_{n+1}(\FF)$. Then $\bx\otimes \xi \perp_f M$ if and only if $\xi M\bx = 0$. 
\end{prop} 

\subsection{The hyperplanes of $\AnF$ which arise from $\veS$}\label{pure tensor} 

For $M\in M_{n+1}(\FF)\setminus\{\lambda I\}_{\lambda\in \FF}$, we put
\begin{equation}\label{eq-id}
\begin{array}{rcl}
\cH_{M, \mathrm{id}_\FF} & := & \veS^{-1}([M^{\perp_f}\cap M_{n+1}^0(\FF)])\\
 & = & \{([\bx], [\xi])\in{\cal P} \mid \xi M \bx = 0\}. 
\end{array}
\end{equation}
(Note that $M^{\perp_f}\cap \Pu(V\otimes V^*) = \{\bx\otimes\xi \neq 0 \mid \xi M\bx = 0\}$, by Proposition \ref{prop1}.) So, $\cH_{M, \idF}$ is a geometric hyperplane of $A_{n,\{1,n\}}(\FF)$. By Proposition \ref{Main0}, it arises from $\veS$. 

In particular, let $\mathrm{rank}(M) = 1$, say $M = \ba\otimes \alpha$ for a pure tensor $\ba\otimes\alpha\in\Pu(V\otimes V^*)$. 
By Proposition \ref{prop0}, if $\cH_{\ba\otimes\alpha, \idF} = \cH_{\bfb\otimes\beta, \idF}$ then $\bfb\otimes\beta \in\langle \ba\otimes\alpha, I\rangle$. However, all rank 1 matrices contained in $\langle \ba\otimes\alpha, I\rangle$ are proportional to $\ba\otimes\alpha$. Hence $\cH_{\ba\otimes\alpha, \idF} = \cH_{\bfb\otimes\beta, \idF}$ if and only if $\ba\otimes\alpha$ and $\bfb\otimes\beta$ are proportional, namely they correspond to the same pair $(a, A)$, where $a = [\ba]$ and $A = [\alpha]$ are the point and the hyperplane of $\PG(n,\FF)$ represented by $\ba$ and $\alpha$ respectively. In fact, as proved in \cite[Proposition 1.3]{Pas}, 
\begin{equation}\label{eq-tensor}
\cH_{\ba\otimes\alpha, \idF} ~=~ \cH_{a,A},
\end{equation} 
where $\cH_{a,A}$ is the quasi-singular hyperplane associated to the pair $(a,A)$ as defined in Section \ref{quasi-sing}. 

\subsection{The hyperplanes of $\AnF$ which arise from $\ve_\sigma$}\label{pure tensor bis}  

Let $\sigma\in \Aut(\FF)\setminus\{\idF\}$. Now, for $M\in M_{n+1}(\FF)\setminus\{\lambda I\}_{\lambda\in \FF}$ we have
\[[M^{\perp_f}]\cap \ve_\sigma({\cal P}) ~ = ~ \{[\bx\otimes\xi] \mid \bx \neq 0 \neq \xi, ~\xi(\bx) = 0 \mbox{ and } \xi M\bx^\sigma = 0\}.\]
We put
\begin{equation}\label{eq-sigma}
\begin{array}{rcl}
\cH_{M, \sigma} & := & \ve_\sigma^{-1}([M^{\perp_f}\cap M_{n+1}(\FF)])\\
 & = & \{([\bx], [\xi])\in{\cal P} \mid \xi M \bx^\sigma = 0\}. 
\end{array}
\end{equation}
The set $\cH_{M, \sigma}$ is a geometric hyperplane of $A_{n,\{1,n\}}(\FF)$. By Proposition \ref{Main0}, it arises from $\ve_\sigma$. 

In particular, let $M = \ba\otimes \alpha$ for $\ba\otimes\alpha\in \Pu(V\otimes V^*)$. Put $a = [\ba]$ and $A := [\alpha^{\sigma^{-1}}]$. 

\begin{prop}\label{prop2}
With $a$ and $A$ as above, we have $\cH_{\ba\otimes\alpha, \sigma} = \cH_{a,A}$. 
\end{prop}
{\bf Proof.} By (\ref{f3}) and (\ref{eq-sigma}) a point $([\bx], [\xi])$ of $\AnF$ belongs to $\cH_{\ba\otimes\alpha,\sigma}$ if and only if $\alpha(\bx^\sigma)\xi(\ba) = 0$. The latter is equivalent to $\alpha^{\sigma^{-1}}(\bx)\xi(\ba) = 0$ which in turn characterizes the points of $\cH_{\ba\otimes\alpha^{\sigma^{-1}}, \idF}$, by (\ref{eq-id}). Therefore $\cH_{\ba\otimes\alpha,\sigma} = \cH_{\ba\otimes\alpha^{\sigma^{-1}},\idF}$. However, with $a = [\ba]$ and $A = [\alpha^{\sigma^{-1}}]$, we have $\cH_{\ba\otimes\alpha^{\sigma^{-1}},\idF} = \cH_{a,A}$ by (\ref{eq-tensor}). Hence $\cH_{\ba\otimes\alpha, \sigma} = \cH_{a,A}$, as claimed.  \hfill $\Box$\\

\noindent
{\bf Proof of Theorem \ref{Main1}, claim (1).} Claim (1) of Theorem \ref{Main1} immediately follows from (\ref{eq-tensor}) and Proposition \ref{prop2}. \hfill $\Box$ 

\section{Proof of Theorem \ref{Main1}, claim (2)}\label{Proof}

In view of (\ref{eq-tensor}) of Section \ref{pure tensor} and Proposition \ref{prop2}, claim (2) of Theorem \ref{Main1} amounts to the following: 

\begin{prop}\label{Main3}
Let $\rho$ and $\sigma$ be distinct automorphisms of $\FF$ with $\sigma\neq \idF$ and let $M, N\in M_{n+1}(\FF)\setminus\{O\}$ be such that
\begin{equation}\label{eq-main}
\xi M\bx^\rho = 0 \mbox{ if and only if } \xi N\bx^\sigma = 0
\end{equation}
for every null-traced pure tensor $\bx\otimes\xi \in \Pu((V\otimes V^*)_0)$. Then $\mathrm{rank}(N) = 1$. 
\end{prop} 

Indeed, (\ref{eq-main}) restricted to the pairs $(\bx, \xi)$ such that $\xi(\bx) = 0$ is equivalent to the equality $\cH_{M,\rho} = \cH_{N,\sigma}$. By Proposition \ref{prop2} and since $\sigma\neq \idF$ by assumption, $\mathrm{rank}(N) = 1$ if and only if $\cH_{N,\sigma}$ is quasi-singular. If $\rho\neq \idF$ the same holds true for $M$ and $\cH_{M,\rho}$. In contrast, let $\rho = \idF$. Then $\cH_{M+\lambda I,\rho} = \cH_{M,\rho}$ for every $\lambda\in \FF$. If $\cH_{M, \rho}$ is quasi-singular then $\mathrm{rank}(M+\lambda I) = 1$ for exactly one choice of $\lambda$. 

Condition (\ref{eq-main}) of Proposition \ref{Main3} can be rephrased as follows:
\begin{equation}\label{eq-main-1}
\langle \bx, M\bx^\rho\rangle ~ = ~ \langle \bx, N\bx^\sigma\rangle, \hspace{5 mm} \forall \bx\in V.
\end{equation}

\subsection{Proof of Proposition \ref{Main3} in the case $n \geq 3$}\label{sec nuova} 

\begin{lemma}\label{nuovo 1}
We have $\dim\langle \bx, \by, M\bx^\rho, M\by^\rho\rangle = \dim\langle \bx, \by, N\bx^\sigma, N\by^\sigma\rangle \leq 3$ for every choice of $\bx, \by \in V$.
\end{lemma}
{\bf Proof.} The equality $\dim\langle \bx, \by, M\bx^\rho, M\by^\rho\rangle = \dim\langle \bx, \by, N\bx^\sigma, N\by^\sigma\rangle$ follows from (\ref{eq-main-1}). We shall prove that $\dim\langle \bx, \by, M\bx^\rho, M\by^\rho\rangle \leq 3$ for every choice of $\bx, \by \in V$. By (\ref{eq-main-1}), for every $\bx\in V$ there exist scalars $\lambda_\bx, \mu_\bx\in \FF$ such that 
\begin{equation}\label{eq nuova 1}
N\bx^\sigma ~ = ~ \bx\lambda_\bx + M\bx^\rho\mu_\bx.
\end{equation}
Hence
\begin{equation}\label{eq nuova 2}
\begin{array}{rcl}
N(\bx+\by)^\sigma & = & \bx\lambda_\bx + \by\lambda_\by + M\bx^\rho\mu_\bx + M\by^\rho\mu_\by\\
 & = &  (\bx+\by)\lambda_{\bx+\by} + M(\bx+\by)^\rho\mu_{\bx+\by}.
\end{array}
\end{equation} 
For a contradiction, suppose that $\bx, \by \in V$ exist such that $\dim\langle\bx, \by,M\bx^\rho, M\by^\rho\rangle = 4$. Then (\ref{eq nuova 2}) implies that  $\lambda_\bx = \lambda_{\bx+\by} = \lambda_\by$ and $\mu_\bx = \mu_{\bx+\by} = \mu_\by$.  Hence $\lambda_\bx =\lambda_\by$ and $\mu_\bx = \mu_\by$. Similarly, 
 $\lambda_{\bx t} =\lambda_\by$ and $\mu_{\bx t} = \mu_\by$ for every $t\in \FF^* := \FF\setminus\{0\}$. Therefore $\lambda_{\bx t} = \lambda_\bx$ and $\mu_{\bx t} = \mu_\bx$ for every $t\in \FF^*$. On the other hand,
\[N(\bx t)^\sigma ~ = ~  \bx t\lambda_{\bx t} + M\bx^\rho t^\rho\mu_{\bx t} ~ = ~ \bx\lambda_\bx t^\sigma + M\bx^\rho t^\sigma.\]
Hence $\lambda_{\bx t} =\lambda_\bx t^{\sigma-1}$ and $\mu_{\bx t} = \mu_\bx t^{\sigma-\rho}$. However $\lambda_\bx = \lambda_{\bx t}$ and $\mu_\bx = \mu_{\bx t}$, as we have already proved. It follows that either $\mu_\bx = 0$ or $t^\sigma = t^\rho$ for every $t\in \FF$. The latter contradicts the hypothesis that $\sigma \neq \rho$. Therefore $\mu_\bx = 0$, namely $N\bx^\sigma \in \langle \bx\rangle$. By (\ref{eq-main-1}) this implies that $M\bx^\rho \in \langle\bx\rangle$, which contradicts the hypothesis that $\dim\langle\bx, \by,M\bx^\rho, M\by^\rho\rangle = 4$. \hfill $\Box$ 

\begin{lemma}\label{nuovo 2}
Let $\ba\in V$ be such that $N\ba^\sigma \not \in\langle\ba\rangle$. Then the $2$-supbspce $L_\ba := \langle \ba, N\ba^\sigma\rangle$ of $V$ contains $N\bx^\sigma$ for every $\bx\in V$.
\end{lemma}
{\bf Proof.} Let $\bx\not\in L_\ba$. By Lemma \ref{nuovo 1}, there exists a unique scalar $\lambda_\bx\in \FF$ and a unique vector $v_\bx\in L_\ba$ such that
\begin{equation}\label{eq nuova 3}
N\bx^\sigma ~ = ~  \bx\lambda_\bx + v_\bx.
\end{equation}
As $n \geq 3$, a vector $\by\in V$ also exists such tat $\dim\langle \bx, \by, L_\ba\rangle = 4$. According to (\ref{eq nuova 3}), scalars $\lambda_\by$ and $\lambda_{\bx+\by}$ and vectors $\bv_\by, \bv_{\bx+\by} \in L_\ba$ exist such that
\begin{equation}\label{eq nuova 4}
\begin{array}{rcl}
 N\by^\sigma & = & \by\lambda_\by + \bv_\by,\\
 N(\bx+\by)^\sigma & = &  (\bx+\by)\lambda_{\bx+\by}+ \bv_{\bx+by} ~  =  \\
= ~ N\bx^\sigma + N\by^\sigma & = & \bx\lambda_\bx + \by\lambda_\by + \bv_\bx + \bv_\by.
\end{array}
\end{equation}
The last two equations of (\ref{eq nuova 4}) yield $\lambda_\bx = \lambda_\by$ ($ =\lambda_{\bx+\by}$). Similarly, $\lambda_{\bx t} = \lambda_\by$ for every $t \in \FF^*$. Consequently, $\lambda_{\bx t} = \lambda_\bx$ for every $t\in \FF^*$. On the other hand, $N(\bx t)^\sigma = N\bx^\sigma t^\sigma = \bx\lambda_\bx t^\sigma + \bv_\bx t^\sigma$. Hence $t\lambda_{\bx t} = \lambda_\bx t^\sigma$, namely $\lambda_{\bx t} = \lambda_\bx t^{\sigma-1}$. 
However $\lambda_{\bx t} = \lambda_\bx$, as previously proved. So, if $\lambda_\bx \neq 0$ we obtain that $t^{\sigma-1} = 1$ for every $t\in \FF^*$, namely $\sigma = \idF$, while $\sigma \neq \idF$ by assumption. Therefore $\lambda_\bx = 0$ and (\ref{eq nuova 3}) implies that  $N\bx^\sigma \in L_\ba$.  \hfill $\Box$

\begin{lemma}\label{nuovo 3}
Suppose that $N\ba^\sigma \not\in \langle\ba\rangle$ for at least one vector $\ba\in V\setminus\{0\}$. Then $\mathrm{rank}(N) = 1$.
\end{lemma}
{\bf Proof.} With $\ba$ as in the hypotheses of the lemma,  Lemma \ref{nuovo 2} implies that the image $\mathrm{Im}(f)$ of the semilinear mapping $f:\bx\rightarrow N\bx^\sigma$ is a subspace of $L_\ba = \langle \ba, N\ba^\sigma\rangle$. Hence the kernel $K := \mathrm{Ker}(f)$ of $f$ (which is the $\sigma^{-1}$-image of the kernel of $N$) has dimension $\dim(K) \geq n-1$ ($\geq 2$ since $n \geq 3$ by assumption). Clearly $L_\ba\not\supseteq K$. This is obvious when $n > 3$. When $n = 3$ this claim follows from the fact that $N\ba^\sigma \neq 0$, since $N\ba^\sigma\not\in \langle\ba\rangle$ by assumption. 

Pick a vector ${\bf k} \in K\setminus L_\ba$ and put $\bfb := \ba+{\bf k}$. Then $N\bfb^\sigma = N\ba^\sigma \not\equiv \bfb$. By Lemma \ref{nuovo 2} with $\ba$ replaced by $\bfb$ we obtain that $\mathrm{Im}(f) \subseteq L_\bfb := \langle \bfb, N\bfb^\sigma\rangle = \langle \bfb, N\ba^\sigma\rangle$. Therefore $\mathrm{Im}(f) \subseteq L_\ba\cap L_\bfb = \langle N\ba^\sigma\rangle$. Accordingly, $\mathrm{rank}(N) = 1$. \hfill $\Box$ \\

In order to finish the proof of Proposition \ref{Main3} when $n \geq 3$ the following case remains to consider: $N\ba^\sigma \in \langle \ba\rangle$ for every $\ba\in V$. In this case $N$ is diagonal, say $N = \mathrm{diag}(n_0,n_1,..., n_n)$. We are assuming that for every vector $\bx= (x_i)_{i=0}^n$ of $V$ there exists a scalar $\lambda \in \FF$ such that $n_ix_i^\sigma = x_i\lambda$ for every $i = 0, 1,..., n$. It is easily to see that this can happen only if $N = O$. However $N\neq O$ by assumption. Hence the hypothesis of Lemma \ref{nuovo 3} holds true and $\mathrm{rank}(N) = 1$ by that lemma. \hfill $\Box$

\subsection{The case $n = 2$}\label{sec riciclo} 

Throught this subsection $n = 2$. In this case condition (\ref{eq-main-1}) is equivalent to the following two conditions, where $(\bx, M\bx^\rho, N\bx^\sigma)$ is the $3\times 3$ matrix with $\bx, M\bx^\rho$ and $N\bx^\sigma$ as the columns: 
\begin{equation}\label{eq-main-1.1}
\mathrm{det}(\bx, M\bx^\rho, N\bx^\sigma) ~ = ~ 0, \hspace{3 mm} \forall\bx\in V;
\end{equation}
\begin{equation}\label{eq-main-1.2}
M\bx^\rho\in \langle \bx\rangle \mbox{ if and only if } N\bx^\sigma\in \langle\bx\rangle, \hspace{3 mm} \forall\bx\in V.
\end{equation}
Let $m_{i,j}$ and $n_{i,j}$ be the $(i,j)$-entries of $M$ and $N$ respectively. For $\{i,j,k\} = \{0,1,2\}$, by (\ref{eq-main-1.1}) with $\bx = \be_k$ we obtain that 
\begin{equation}\label{eq-main-1.3}
m_{i,k}n_{j,k} ~ = ~ m_{j,k}n_{i,k}. 
\end{equation}
Equality (\ref{eq-main-1.3}) shows that there exist a matrix $R = (r_{i,j})_{i,j=0}^2\in M_3(\FF)$ and diagonal matrices 
\[\begin{array}{cc}
\Delta_M = \mathrm{diag}(\mu_0, \mu_1, \mu_2), & \Delta_N = \mathrm{diag}(\nu_0, \nu_1, \nu_2),\\
D_M =  \mathrm{diag}(m_0, m_1, m_2), & D_N = \mathrm{diag}(n_0, n_1, n_2)
\end{array}\]
such that $r_{i,i} = 0$ for every $i = 0,1,2$ and 
\begin{equation}\label{eq-main-2}
M ~=~ R\Delta_M+D_M, \hspace{5 mm} N ~ = ~ R\Delta_N + D_N.
\end{equation}
So, if ${\bf r}_k = (r_{i,k})_{i=0}^2$ is the $(k+1)$-th column of $R$ and $\bm_k$, $\bn_k$ are the $(k+1)$-th columns of $M$ and $N$ respectively, we have $\bm_k = \br_k\mu_k$ (namely $m_{i,k} = r_{i,k}\mu_k$ for every $i \neq k$), $\bn_k = \br_k\nu_k$ (namely $n_{i,k} = r_{i,k}\nu_k$ for every $i\neq k$), $m_k = m_{k,k}$ and $n_k = n_{k,k}$ for every $k = 0, 1, 2$. Moreover (\ref{eq-main-1.2}) implies that we can choose $R$, $\Delta_M$ and $\Delta_N$ in such a way that
\begin{equation}\label{eq-main-3}
\br_k \neq 0 \mbox{ if and only if } \mu_k \neq 0 \mbox{ if and only if } \nu_k\neq 0.
\end{equation}
A matrix $R$ as above will be called a {\em skeleton} of the pair $(M,N)$. Henceforth $R$ will always stand for a given skeleton of $(M,N)$. We warn that if $R \neq O$ then $R$ is not uniquely determined by $M$ and $N$. Indeed $M$ and $N$ determine $R$ up to right multiplication $RD$ with $D = \mathrm{diag}(\lambda_0, \lambda_1, \lambda_2)$ such that, for $k \in \{0,1,2\}$, if $\br_k \neq 0$ then $\lambda_k \neq 0$.       

Suppose firstly that $R = O$. Then $M = D_M$ and $N = D_N$. With the help of the Identity Principle for semi-polynomials (Theorem \ref{Princ}) we see that $N \equiv \be_k\otimes\eta_k$ for some $k \in \{0,1,2\}$. Hence $\mathrm{rank}(N) = 1$. For the rest of this subsection we assume that $R\neq O$. 

\subsubsection{The case where $R \neq O$ but two columns of $R$ are null}

Suppose $R\neq O$ but two of its columns are null. To fix ideas, let $\br_0 \neq 0 = \br_1 = \br_2$. If $\br_0 \neq \be_1$, we can replace the natural basis $E = (\be_0, \be_1, \be_2)$ of $V$ with a new basis $E' = (\be'_0, \be'_1, \be'_2)$ such that $\be'_0 = \be_0$, $\be'_1 = \br_0$ and $\langle \be'_1, \be'_2\rangle = \langle \be_1, \be_2\rangle$. (Recall that $\br_0 \in \langle \be_1, \be_2\rangle$, since $r_{0,0} = 0$). If $C\in M_3(\FF)$ is the matrix of this change of bases, then we must replace $M$ and $N$ in (\ref{eq-main-1}) with $C^{-1}MC^\rho$ and $C^{-1}NC^\sigma$ respectively. However $C^{-1}NC^\sigma$ and $N$ have the same rank. Hence there is no loss in assuming that $\br_0 = \be_1$.  

So, let $\br_0 = \be_1$. Accordingly, $R = E_{1,0} = \be_1\otimes\eta_0$. Suppose firstly that $\rho = \idF$. From (\ref{eq-main-1.1}) and the Identity Principle (Theorem \ref{Princ}) we get 
\begin{equation}\label{eq-main-4}
\begin{array}{l}
(m_i-m_0)n_j = 0 \mbox{ for } \{i,j\} = \{1,2\}, \\
 (m_i-m_j)n_0 = 0 \mbox{ for } \{i, j\} = \{1,2\},\\  
\mu_0n_2 = 0, \\
 (m_2-m_0)\nu_0 + \mu_0n_0 = 0. 
\end{array}
\end{equation} 
However $\mu_0\neq 0$ by (\ref{eq-main-3}), since $\br_0\neq 0$. So, the third condition of (\ref{eq-main-4}) implies that $n_2 = 0$. Suppose that $n_1 \neq 0$. Then $m_2 = m_0$ by the first of (\ref{eq-main-4}). The fourth of (\ref{eq-main-4}) now yields $\mu_0n_0 =0$, hence $n_0 = 0$ because $\mu_0\neq 0$ (by (\ref{eq-main-3}) and since $\br_0 \neq 0$). Therefore at least one of $n_0$ or $n_1$ is equal to $0$. Consequently $\mathrm{rank}(N) = 1$.  

Let $\rho\neq \idF$. From (\ref{eq-main-1.1}) and the Identity Principle now we get $n_2 = m_2 = 0$, $m_0\nu_0 = n_0\mu_0$ and $m_0n_1 = m_1n_0 = 0$. These relations imply that $\mathrm{rank}(M) = \mathrm{rank}(N) = 1$. 

\subsubsection{The case where at most one column of $R$ is null} 

Suppose that at most one column of $R$ is null. To fix ideas, let $\br_i \neq 0$ for $i \leq s$ where $s \in \{1,2\}$ and $\br_2 = 0$ if $s = 1$. Therefore $\mu_i \neq 0 = \mu_j$ and $\nu_i \neq 0 = \nu_j$ for $i \leq s < j$. Accordingly, $\bn_i \neq 0$ for every choice of $i\leq s$. 

\begin{lemma}\label{nuovo}
If $\bn_i \equiv \bn_j$ for every choice of $i < j \leq s$ then $\mathrm{rank}(N)= 1$. 
\end{lemma}
{\bf Proof.} Suppose that $\bn_i \equiv \bn_j$ for every choice of $i < j \leq s$. If $s = 2$ then clearly $\mathrm{rank}(N) = 1$. Suppose $s = 1$. By (\ref{eq-main-1.1}) via the Identity Principle we obtain the following among a number of other relations: either $n_2 = 0$ or $r_{0,1} = r_{0,1} = 0$. If $n_2 = 0$ then again $\mathrm{rank}(N) = 1$. 

Suppose that $n_2 \neq 0$. Then $r_{0,1} = r_{1,0} = 0$. Accordingly, $n_{0,1} = n_{1,0} = 0$ and therefore $n_0 = n_1 = 0$ since $\bn_0 \equiv \bn_1$. Consequently, the first two rows of $N$ are null. Hence $\mathrm{rank}(N) = 1$.   \hfill $\Box$  

\begin{lemma}\label{prov1}
We have $\bn_0 \not \equiv \bn_1$ if and only if $\br_0 = \be_1$, $\br_1 = \be_0$ and $n_0n_1 \neq \nu_0\nu_1$. 
\end{lemma}
{\bf Proof.} Suppose that $\rho = \idF$. By (\ref{eq-main-1.1}) with $\bx\in \langle \be_0, \be_1\rangle$ we obtain the following relations:
\begin{equation}\label{eq-two-1}
\begin{array}{rcl}
r_{2,1}r_{1,0}\nu_1 & = & r_{2,0}n_1,\\
r_{2,0}(\nu_0(m_1-m_0) + \mu_0n_0) & = & r_{2,1}r_{1,0}\mu_1\nu_0,\\
r_{2,1}(\nu_1(m_0-m_1) + \mu_1n_1) & = & r_{2,0}r_{0,1}\mu_0\nu_1,\\
r_{2,0}r_{0,1}\nu_0 & = & r_{2,1}n_0.
\end{array}
\end{equation} 
Suppose that $r_{1,0}\neq 0$. Then $n_{2,1} = r_{2,1}\nu_1 = r_{2,0}n_1/r_{1,0} = r_{2,0}\nu_0n_1/r_{1,0}\nu_0 = n_{2,0}\lambda_{0,1}$, where $\lambda_{0,1} := n_1/n_{1,0}$. (The equality $r_{2,1}\nu_1 = r_{2,0}n_1/r_{1,0}$ follows from the first equation of (\ref{eq-two-1}).) Trivially, $n_{1,1} = n_1 = n_{1,0}\lambda_{0,1}$. Therefore 
\begin{equation}\label{eq-two-2.1}
n_{i,1} = n_{i,0}\lambda_{0,1} \mbox{ for } i \in \{1,2\}, \mbox{ where } \lambda_{0,1} = \frac{n_1}{n_{1,0}}.
\end{equation} 
Similarly, if $r_{0,1} \neq 0$ then by exploiting the last equation of (\ref{eq-two-1}) we obtain: 
\begin{equation}\label{eq-two-2.2}
n_{i,0} = n_{i,1}\lambda_{1,0} \mbox{ for } i \in \{0,2\}, \mbox{ where } \lambda_{1,0} = \frac{n_0}{n_{0,1}}.
\end{equation}
By (\ref{eq-two-2.1}), if $n_{2,1}\neq 0$ then $\lambda_{0,1}\neq 0$ and $n_{2,0}\neq 0$. Similarly, by (\ref{eq-two-2.2}), if $n_{2,0}\neq 0$ then $\lambda_{1,0} \neq 0 \neq n_{2,1}$. So, either $\lambda_{0,1} \neq 0 \neq \lambda_{1,0}$ and $n_{2,1} \neq 0 \neq n_{2,0}$, or $n_{2,0} = n_{2,1} = 0$, hence $r_{2,0} = r_{2,1} = 0$. In the latter case, up to replace $\nu_0, \nu_1, \mu_0$ and $\mu_1$ with suitable non-zero scalars, we can assume that $\br_0 = \be_1$ and $\br_1 = \be_0$. 

Suppose the first case occurs: $n_{2,0} \neq 0 \neq n_{2.1}$. Then $\lambda_{1,0} \neq 0 \neq \lambda_{0,1}$ and $n_{2,1} = n_{2,0}\lambda_{0,1} = n_{2,1}\lambda_{1,0}\lambda_{0,1}$. Therefore $\lambda_{1,0}\lambda_{0,1} = 1$. Hence $n_{1,1} = n_{1,0}\cdot n_{0,1}/n_0 = n_{1,0}\lambda_{1,0}^{-1} = n_{1,0}\lambda_{0,1}$. So, $n_{i,1} = n_{i,0}\lambda_{0,1}$ for every $i$,  namely $\bn_1 = \bn_0\lambda_{0,1}$. In this case $\bn_1 \equiv \bn_0$. 

In the second case, where $\br_0 = \be_1$ and $\br_1 = \be_0$, we have $\bn_1 \equiv \bn_0$ if and only if $n_0n_1 = n_{0,1}n_{1,0}$ ($= \nu_1\nu_0$, since $r_{0,1} = r_{1,0} = 1$), as claimed in the statement of the lemma. 

Suppose now that at least one of $r_{1,0}$ or $r_{0,1}$ is null. Let $r_{1,0} = 0$, to fix ideas. Then $r_{2,0}\neq 0$, since $\br_0 \neq 0$ by assumption. The first and second equation of (\ref{eq-two-1}) now yield $n_1 = 0$ and $m_0-m_1 = \mu_0n_0/\nu_0$ respectively. As  $n_1 = 0$, the third equation of (\ref{eq-two-1}) yields 
\begin{equation}\label{sequel 1}
r_{2,1}(m_0-m_1) = r_{2,0}r_{0,1}\mu_0.
\end{equation} 
Suppose that $r_{0,1} \neq 0$. Then $m_0 \neq m_1$, otherwise (\ref{sequel 1}) forces $r_{2,0} = 0$ while $r_{2,0}\neq 0$. Hence $r_{1,1} = r_{2,0}r_{0,1}\mu_0/(m_0-m_1)$. However $m_0-m_1 = \mu_0n_0/\nu_0$. Hence $n_0 \neq 0$ and $r_{i,1} = r_{i,0}r_{0,1}\nu_0/n_0$, namely $n_{i,1} = n_{i,0}\cdot n_{0,1}/n_0$. This holds for $i = 2$ as well as for $i = 0$, as we already know from (\ref{eq-two-2.2}). However this relation also holds for $i = 1$, as $n_{1,1} = n_1 = 0$ and $n_{1,0} = r_{1,0}\nu_0 = 0$, since $r_{1,0} = 0$ by assumption. So, $\bn_1\equiv \bn_0$.  

On the other hand, still assuming that $r_{1,0} = 0$, suppose that $r_{0,1} = 0$. Then $r_{2,1} \neq 0$ since $\br_1 \neq 0$ and, by the fourth equation of (\ref{eq-two-1}), we also have $n_0 = 0$. So, $\br_0 = \br_1 = \be_2$, $\bn_0 = \be_2\nu_0$ and $\bn_1 = \be_2\nu_1$. Again, $\bn_0 \equiv \bn_1$. 

So far we have proved that the statement of the lemma holds true when $\rho = \idF$. Let $\rho \neq \idF$.  By (\ref{eq-main-1.1}) with $\bx\in \langle \be_0, \be_1\rangle$ we now obtain the following relations: 
\begin{equation}\label{eq-two-1 bis}
\begin{array}{ll}
r_{2,1}r_{1,0}\nu_1 =  r_{2,0}n_1, & r_{2,1}r_{1,0}\mu_1 = r_{2,0}m_1,\\
r_{2,0}r_{0,1}\nu_0  =  r_{2,1}n_0, & r_{2,0}r_{0,1}\mu_0  =  r_{2,1}m_0, \\
r_{2,0}\nu_0m_0  =  r_{2,0}\mu_0n_0, & r_{2,1}\nu_1m_1  =  r_{2,1}\mu_1n_1. 
\end{array}
\end{equation} 
As when $\rho = \idF$, if $r_{1,0}\neq 0 \neq r_{0,1}$ and either $r_{2,0} \neq 0$ or $r_{2,1} \neq 0$ then $\bn_1\equiv \bn_0$ and $\bm_1\equiv \bm_0$. If $r_{1,0} = 0 = r_{0,1}$ we can assume that $\br_0 = \be_1$ and $\br_1 = \be_0$. In this case $\bn_1\equiv \bn_0$ if and only if $\nu_0\nu_1 = n_0n_1$. Similarly, with $\br_0 = \be_1$ and $\br_1 = \be_0$ we have $\bm_1 \equiv \bm_0$ if and only if $m_0m_1 = \mu_0\mu_1$. 

Suppose that at least one of $r_{1,0}$ or $r_{0,1}$ is null. Let $r_{1,0} = 0$, to fix ideas. Then $r_{2,0}\neq 0$ since $\br_0 \neq 0$. The first and second equation of (\ref{eq-two-1 bis}) now yield $n_1 = m_1 = 0$ while the fifth equation yields $\nu_0m_0 = \mu_0n_0$. Suppose that $r_{0,1} \neq 0$. Then the third and fourth equation yield $n_0 \neq 0 \neq m_0$. Therefore $n_{i,1} = n_{i,0}\cdot n_{0,1}/n_0$ and $m_{i,1} = m_{i,0}\cdot m_{0,1}/m_0$ for every $i = 0, 1, 2$. When $i = 2$ these equalities follow from the first two relations of (\ref{eq-two-1 bis}), they are trivial for $i = 0$ and when $i = 1$ they hold because $r_{1,0} = n_1 = m_1 = 0$. Hence $\bn_1 \equiv \bn_0$ and $\bm_1 \equiv \bm_0$. Similarly, if $r_{0,1} = 0 \neq r_{1,0}$ then $\bn_1 \equiv \bn_0$ and $\bm_1 \equiv \bm_0$.

Finally, let $r_{1,0} = r_{0,1} = 0$. Then $r_{2,0} \neq 0 \neq r_{2,1}$ because $\br_0 \neq 0 \neq \br_1$. In this case (\ref{eq-two-1 bis}) force $n_0 = m_0 = n_1 = m_0 = m_1 = 0$. Hence $\bn_0\equiv \bn_1$.   \hfill $\Box$   

\begin{lemma}\label{prov2}
Let $\br_0 = \be_1$, $\br_1 = \be_0$ and $\br_2 = 0$. Then $\bn_0 \equiv \bn_1$. 
\end{lemma}
{\bf Proof.} Suppose firstly that $\rho = \idF$. By (\ref{eq-main-1.1}) we obtain that $n_2 = 0$ and 
\begin{equation}\label{sequel 2}
\begin{array}{rclcrcl}
\nu_0(m_0-m_2) & = & \mu_0n_0, & & \nu_1(m_1-m_2) & = & \mu_1n_1.\\
n_0(m_1-m_2) & = & \mu_1\nu_0, & & n_1(m_0-m_2) & = & \mu_0\nu_1.
\end{array}
\end{equation}
The third and fourth equation of (\ref{sequel 2}) imply $n_1 \neq 0 \neq n_0$. From the first and fourth equation (the second and the third) we obtain $m_0-m_2  = \mu_0n_0/\nu_0 = \mu_0\nu_1/n_1$ (respectively $m_1-m_2  = \mu_1n_1/\nu_1 = \mu_1\nu_0/n_0$). Hence $n_0/\nu_0 = \nu_1/n_1$, namely $n_0n_1 = \nu_0\nu_1$. Therefore $\bn_0 \equiv \bn_1$, by Lemma \ref{prov1}. 

Let now $\rho \neq \idF$. By (\ref{eq-main-1.1}) we obtain that $m_2 = n_2 = 0$ and
\begin{equation}\label{sequel 2 bis}
\begin{array}{cc}
\nu_0m_0  =  \mu_0n_0, & \nu_1m_1  = \mu_1n_1,\\
m_0n_1  =  \mu_0\nu_1, & m_1n_0  =  \mu_1\nu_0.
\end{array}
\end{equation}
The first and second equations of (\ref{sequel 2 bis}) yield $n_0 = m_0\nu_0/\mu_0$ and $n_1 = m_1\nu_1/\mu_1$. From these equalities and the third and fourth equations of (\ref{sequel 2 bis}) we get $m_0m_1 = \mu_0\mu_1$. Similarly, $n_0n_1 = \nu_0\nu_1$. Hence $\bm_0 \equiv \bm_1$ and $\bn_0 \equiv \bn_1$ by Lemma \ref{prov1}. \hfill $\Box$

\begin{lemma}\label{prov3}
We have $\bn_0 \equiv \bn_1$. 
\end{lemma}
{\bf Proof.} For a contradiction, suppose that $\bn_0\not\equiv \bn_1$. Then $\br_0 = \be_1$ and $\br_1 = \be_0$ by Lemma \ref{prov1}. The column $\bn_2$, which is non-null by Lemma \ref{prov2}, cannot be proportional to both $\bn_0$ and $\bn_1$. Therefore, by Lemma \ref{prov1}, either $\br_0 = \be_2$ and $\br_2 = \be_0$ or $\br_1 = \be_2$ and $\br_2 = \be_1$. Each of these two cases crashes against the fact that $\br_0 = \be_1$ and $\br_1 = \be_0$. \hfill $\Box$ \\

If $s = 1$ the conclusion $\mathrm{rank}(N) = 1$ follows from Lemmas \ref{nuovo} and \ref{prov3}. When $s = 2$ the same argument used to conclude that $\bn_0 \equiv \bn_1$ also yields $\bn_1\equiv \bn_2$. Again, $\mathrm{rank}(N) = 1$ by Lemma \ref{nuovo}. \hfill $\Box$  

\section{Proof of Theorem \ref{Main4}}\label{Absolute} 

In the proof of the next lemma $\perp$ is the collinearity relation of $\AnF$ and, for a point $(x, X)$ of $\AnF$, $(x,X)^\perp = {\cal M}_x\cup{\cal M}_X$ is the set of points of $\AnF$ collinear with $(x,X)$ or equal to it. Recall that two points $(x,X)$ and $(y,X)$ of $\AnF$ are collinear if and only if either $x = y$ or $X = Y$; they are at distance at most 2 if and only if either $x\in Y$ or $y\in X$ (possibly both $x\in Y$ and $y\in X$).  

\begin{lemma}\label{V0}
Let $(a,A)$ and $(b,B)$ be points of $\AnF$ at distance $3$ in the collinearity graph of $\AnF$. Then $\cH_{a,A}\cap\cH_{b,B}$, which is a hyperplane of the singular hyperplane $\cH_{a,A}$ of $\AnF$, is a maximal subspace of $\cH_{a,A}$.
\end{lemma}
{\bf Proof.} Put $\cH := \cH_{a,A}$, ${\cal I} := \cH_{a,A}\cap\cH_{b,B}$ and ${\cal C} := \cH\setminus{\cal I}$. We shall prove that $\perp$ induces a connected graph on $\cal C$. The conclusion will follow from Shult \cite[Lemma 4.1.1]{S}. 

Recall that, for a point $(p,H)$ of $\AnF$, the hyperplane $\cH_{p,H}$ is the set of points at distance at most $2$ from $(p,H)$. Accordingly, $\cH$ is the set of points at distance at most 2 from $(a,A)$ and $\cal I$ is the set of points at distance at most $2$ from both $(a,A)$ and $(b,B)$. Hence $(a,A)\in {\cal C}$, as $(b,B)$ has distance 3 from $(a,A)$ by assumption. We shall prove that every point $(c,C)\in{\cal C}$ can be connected to $(a,A)$ by a path fully contained in $\cal C$. 

If either $(c,C)\perp (a,A)$ or $(c,C)^\perp\cap (a,A)^\perp\not\subseteq{\cal I}$ there is nothing to prove. Suppose that $(c,C)$ is at distance 2 from $(a,A)$. Then either $c\in A$ or $a\in C$ (possibly $c\in A$ and $a\in C$). Let $c\in A$, to fix ideas. Then $(c,A)\in (a,A)^\perp\cap(c,C)^\perp$. As both $(a,A)$ and $(c,C)$ are at distance 3 from $(b,B)$ (recall that $(c,C)\in {\cal C}$), neither $b\in A$ nor $c\in B$. Accordingly, $(c,A)$ is at distance 3 from $(b,B)$, namely $(c,A) \in {\cal C}$. Hence $(c,C)\perp(c,A)\perp(a,A)$ is a path in $\cal C$ from $(c,C)$ to $(a,A)$.  \hfill $\Box$ \\

Throughout the rest of this section we assume that $|\Aut(\FF)| > 1$. For $\sigma\in\Aut(\FF)$, if $\cH$ is a hyperplane of $\AnF$ which arises from $\sigma$, we denote by $\ve_{\sigma \mid \cH}$ the embedding of $\cH$ in $\langle \ve_\sigma(\cH)\rangle$ induced by $\ve_\sigma$ on $\cH$. 

Recall that, by Theorem \ref{Main1}, the quasi-singular hyperplanes of $\AnF$ arise from $\ve_\sigma$ for every $\sigma \in \Aut(\FF)$ while the hyperplanes which are not quasi-singular arise from $\ve_\sigma$ for at most one $\sigma\in \Aut(\FF)$. Recall also that $\veS = \ve_\sigma$ with $\sigma = \idF$.  

\begin{lemma}\label{V1}
For every singular hyperplane $\cH_1$ of $\AnF$ there exists a hyperlane $\cH_2$ of $\AnF$ such that:
\begin{itemize}
\item[$(1)$] the intersection $\cH_1\cap \cH_2$ is maximal as a subspace of $\cH_1$ and, regarded as a hyperlane of $\cH_1$, it arises from $\ve_{\sigma \mid \cH_1}$ for every $\sigma\in \Aut(\FF)$;
\item[$(2)$] the hyperplane $\cH_2$ arises from $\ve_\rho$ for a unique $\rho\in\Aut(\FF)$ and $\cH_1\cap \cH_2$, regarded as a hyperplane of $\cH_2$,  arises from $\ve_{\rho\mid\cH_2}$.
\end{itemize}   
\end{lemma}
{\bf Proof.} Given two points $(a,A)$ and $(b,B)$ of $\AnF$ ad distance 3, let $\cH_1 = \cH_{a,A}$ and $\cH_1' = \cH_{b,B}$. Then $\cH_1$ and $\cH'_1$ arise from $\ve_\sigma$ for every $\sigma\in\Aut(\FF)$. Moreover $\cH_1\cap\cH'_1 = \ve^{-1}_\sigma(\langle \ve_\sigma(\cH_1)\rangle\cap\langle\ve_\sigma(\cH'_1)\rangle)$. By Lemma \ref{V0} the subspace ${\cal K}:= \cH_1\cap\cH'_1$ is a hyperplane as well as maximal subspace of both $\cH_1$ and $\cH'_1$. Hence $\cal K$ arises from both $\ve_{\sigma \mid \cH_1}$ and $\ve_{\sigma\mid \cH'_1}$ by Proposition \ref{Gen 1}. Explictly, for every $\sigma\in \Aut(\FF)$ we have
\begin{equation}\label{eqV1}
\langle\ve_\sigma({\cal K})\rangle ~ = ~ \langle \ve_\sigma(\cH_1)\rangle\cap\langle\ve_\sigma(\cH'_1)\rangle \mbox{ and } {\cal K} = \ve_\sigma^{-1}(\langle\ve_\sigma({\cal K})\rangle). 
\end{equation}
Let $\ba$ and $\bfb$ be representative vectors of $a$ and $b$ respectively and $\alpha$ and $\beta$ linear functional describing $A$ and $B$ and choose an automorphism $\rho\in\Aut(\FF)$. We know from Section \ref{pure tensor} and Proposition \ref{prop2} that $\cH_1 = \cH_{\ba\otimes\alpha^\rho, \rho}$ and $\cH_1' = \cH_{\bfb\otimes \beta^\rho, \rho}$. As we have chosen $(a,A)$ and $(b,B)$ at mutual distance 3, neither $\ba\equiv\bfb$ nor $\alpha\equiv\beta$. Hence all matrices of $\langle \ba\otimes\alpha^\rho, \bfb\otimes\beta^\rho\rangle$ have rank 2 except those which are proportional to either $\ba\otimes\alpha^\rho$ or $\bfb\otimes\beta^\rho$. Let $M$ be one of the matrices of $\langle \ba\otimes\alpha^\rho, \bfb\otimes\beta^\rho\rangle$ of rank $2$ and put $\cH_2 := \cH_{M, \rho}$ (notation as in Section \ref{pure tensor bis}). Then $\cH_2$ is not quasi-singular. Hence it arises from $\ve_\rho$ but it does not arise from $\ve_\sigma$ for any $\sigma\neq \rho$, by Theorem \ref{Main1}. Moreover $\cH_2\cap \cH_1 = \cH_2\cap \cH'_1 = \cH_1\cap \cH'_1 = {\cal K}$. These equalities combined with (\ref{eqV1}) show that the hyperplane $\cH_1\cap\cH_2$ of $\cH_2$ arises from  $\ve_{\rho\mid\cH_2}$.   \hfill $\Box$

\begin{lemma}\label{V2}
Given $\cH_1$, $\cH_2$ and $\rho$ as in the statement of Lemma \ref{V1} and chosen an automorphism $\sigma$ of $\FF$ different from $\rho$, there exists a hyerplane $\cH_3$ of $\AnF$ which arises from $\ve_\sigma$ but not from $\ve_\rho$ and such that $\cH_3\cap\cH_2$ properly contains $\cH_1\cap \cH_2$. Moreover, $\cH_1\cap \cH_3 = \cH_1\cap\cH_2$ and $\ve_\sigma(\cH_1\cap\cH_2)$ spans $\langle\ve_\sigma(\cH_1)\rangle\cap\langle\ve_\sigma(\cH_3)\rangle$  while $\ve_\sigma(\cH_3\cap\cH_2)$ spans $\langle \ve_\sigma(\cH_3)\rangle$.  
\end{lemma}
{\bf Proof.} Let $\Sigma_\sigma$ be the codomain of $\ve_\sigma$, namely $\Sigma_\sigma = \PG(M_{n+1}^0(\FF))$ if $\sigma = \idF$ and $\Sigma_2 = \PG(M_{n+1}(\FF))$ if otherwise.

Let $\cal X$ be the collection of subspaces $X$ of $\cH_2$ which contain $\cH_1\cap \cH_2$ and such that $\ve_\sigma(X)$ does not span $\Sigma_\sigma$. Note that $\cH_1\cap \cH_2\in {\cal X}$. By Zorn's Lemma we can see that every member of $\cal X$ is contained in maximal member of $\cal X$. Let $X\in {\cal X}$ be one of them. Then $X \subset \cH_2$. Indeed if otherwise then $\ve_\sigma(\cH_2)$ is contained in a hyperplane $H$ of $\Sigma_\sigma$. Hence $\ve_\sigma^{-1}(H)$ is a hyperplane of $\AnF$ and contains $\cH_2$. However the hyperplanes of $\AnF$ are maximal subspaces. Therefore $\cH_2 = \ve_\sigma^{-1}(H)$. Consequently, $\cH_2$ arises from both $\ve_\rho$ and $\ve_\sigma$. This contradicts condition (2) of Lemma \ref{V1}. 

So, $X\subset \cH_2$. Pick $x\in \cH_2\setminus X$. Then $\ve_\sigma(X\cup\{x\})$ spans $\Sigma_\sigma$, by the maximality of $X$ in $\cal X$. This shows that $H := \langle\ve_\sigma(X)\rangle$ is a hyperplane of $\Sigma_\sigma$. Put $\cH_3 := \ve_\sigma^{-1}(H)$. Then $\cH_3$ is a hyperplane of $\AnF$ and arises from $\ve_\sigma$. As $\cH_1$ does not arise from $\ve_\sigma$, we have $\cH_3\neq \cH_1$. The intersection $\cH_3\cap \cH_2$ contains $X$. It cannot be larger than $X$, otherwise $\ve_\sigma(\cH_3\cap \cH_2)$ spans $\Sigma_\sigma$ by the maximality of $X$, while $\ve_\sigma(\cH_3\cap\cH_2)$ is contained in the hyperplane $H$ of $\Sigma_\sigma$, by defnition of $\cH_3$. Therefore $\cH_3\cap\cH_2 = X$. Clearly $X \supset \cH_1\cap \cH_2$. Indeed $\ve_\sigma(\cH_1\cap\cH_2)$ spans a huperplane of the hyperplane $\langle \ve_\sigma(\cH_1)\rangle$ of $\Sigma_\sigma$ while $\ve_\sigma(X)$ spans the hyperplane $H = \langle \ve_\sigma(\cH_3)\rangle$ of $\Sigma_2$. 

We have $\langle\ve_\sigma(\cH_3\cap\cH_2)\rangle = \langle\ve_\sigma(\cH_3)\rangle$ by definiton of $\cH_3$. On the other hand, $\langle\ve_\sigma(\cH_3\cap\cH_1)\rangle$ contains $\langle\ve_\sigma(\cH_1\cap\cH_2)\rangle$, which is a hyperplane of $\langle \ve_\sigma(\cH_1)\rangle$ by (1) of Lemma \ref{V1}. As $\cH_1\cap \cH_2\subset \cH_3$, the subspace $\langle \ve_\sigma(\cH_1\cap\cH_2)\rangle$ is also a hyperplane of $\langle\ve_\sigma(\cH_3)\rangle$ and must be the same as $\langle\ve_\sigma(\cH_1)\rangle\cap\langle\ve_\sigma(\cH_3)\rangle$. Consequently, $\cH_1\cap\cH_2 = \cH_1\cap \cH_3$.  

Finally, the hyperplane $\cH_3$ does not arise from $\ve_\rho$. Indeed $\ve_\rho(\cH_1\cap \cH_2)$ spans a hyperplane of $\langle \ve_\rho(\cH_2)\rangle$ by condition (2) of Lemma \ref{V1} while $\cH_3$ properly contains $\cH_1\cap \cH_2$. So, if $\cH_3$ also arises from $\ve_\rho$ then $\cH_3 \supseteq \cH_2$, while $\cH_3\cap \cH_2 = X \subset \cH_2$.   \hfill $\Box$ \\

\noindent
{\bf End of the proof of Theorem \ref{Main4}.} For a contradiction, suppose that $\AnF$ admits an embedding $\vet:\AnF\rightarrow\widetilde{\Sigma}$ which covers both $\ve_\sigma$ and $\ve_\rho$ for two distinct automorphisms $\sigma$ and $\rho$ of $\FF$. If a hyperplane $\cH$ arises from $\ve_\sigma$ or $\ve_\rho$ then it also arises from $\vet$, by Proposition \ref{Gen 2} and since all hyperplanes of $\AnF$ are maximal subspaces \cite[Theorem 1.5]{Pas}. The restriction $\vet_{\mid\cH}$ of $\vet$ to $\cH$ covers $\ve_{\sigma\mid\cH}$. Consequently, if $\cal K$ is a hyperplane of $\cH$ which arises from $\ve_{\sigma\mid\cH}$ and is maximal as a subspace of $\cH$, then $\cal K$ also arises from $\vet_{\mid\cH}$, by Proposition \ref{Gen 2}. Of course, the same claims hold true for the embedding $\ve_\rho$ and the hyperplanes which arise from it. 

Let now $\cH_1, \cH_2$ and $\cH_3$ be as in Lemmas \ref{V1} and \ref{V2}. The singular hyperplane $\cH_1$ arises from both $\ve_\sigma$ and $\ve_\rho$, the hyperlane $\cH_2$ arises from $\ve_\rho$ and $\cH_3$ arises from $\ve_\sigma$. Hence each of these hyperplanes arises from $\vet$. Moreover, the hyperplane ${\cal K} := \cH_1\cap\cH_2$ of $\cH_1$ arises from $\vet_{|\cH_1}$, by condition (1) of Lemma \ref{V1}. It follows that $K := \langle\vet({\cal K})\rangle$ is a hyperplane of the hyperplane $H_1 := \langle \vet(\cH_1)\rangle$ of $\widetilde{\Sigma}$. Moreover $K$ is also contained in the hyperplanes $H_2 := \langle \vet(\cH_2)\rangle$ and $H_3 := \langle\vet(\cH_3)\rangle$. Clearly $K$ is a hyperplane of both $H_2$ and $H_3$. Moreover, $H_2\neq H_3$, since $\cH_2\neq\cH_3$. Therefore $H_2\cap H_3 = K$. However this conclusion crashes against the fact that, as stated in Lemma \ref{V2}, the intersection $\cH_1\cap \cH_3$ properly contains ${\cal K} = \cH_1\cap \cH_3$. Therefore no embedding of $\AnF$ covers both $\ve_\sigma$ and $\ve_\rho$.  \hfill $\Box$

\section{Quotients of $\veS$ and its twistings}\label{sec last}

Consider $\ve_\sigma:\AnF\rightarrow\PG(W)$ with $\sigma\in\Aut(\FF)$, where $W = M_{n+1}^0(\FF)$ if $\sigma = \idF$ and  $W = M_{n+1}(\FF)$ if otherwise. 

\begin{lemma}\label{quot1}
Let $S$ be a subspace of $W$. If all matrices of $S\setminus\{O\}$ have rank at least $3$ then the corresponding subspace $[S]$ of $\PG(W)$ defines a quotient of $\ve_\sigma$. When $\sigma = \idF$ the converse also holds: $[S]$ defines a quotient of $\veS$ only if all matrices of $S\setminus\{O\}$ have rank at least $3$.   
\end{lemma}
{\bf Proof.} The $\ve_\sigma$-images of the points of $\AnF$ are represented by matrices of rank 1. The sum of two matrices of rank 1 has rank at most 2. 
The first claim of the lemma follows from this remark and condition (\ref{eq quot}) of Section \ref{prel3 sub}. Turning to the second claim, every null-traced matrix of rank 1 represents the $\veS$-image of a point of $\AnF$. The second claim of the lemma follows from the well known fact that every matrix of rank 2 is the sum of two matrices of rank 1. \hfill $\Box$\\

By the first part of Lemma \ref{quot1}, for every non-singular matrix $M\in W$, the point $[M]$ of $\PG(W)$ defines a quotient of $\ve_\sigma$. In particular, $\ve_\sigma/[I]$ is defined for every $\sigma\neq\idF$ while $\veS/[I]$ is defined if and only if $\mathrm{char}(\FF)$ is positive and divides $n+1$. Indeed $\mathrm{trace}(I) = 0$ if and only if $\mathrm{char}(\FF)$ is positive and divides $n+1$. 

We say that an embedding of $\AnF$ is {\em polarized} if all singular hyperplanes of $\AnF$ arise from it. For instance, as we know from Theorem \ref{Main1}, the natural embedding of $\AnF$ and all of its twistings are polarized.   

\begin{note}
\em
The choice of the word ``polarized" is suggested by the theory of dual polar spaces, where a projective embedding of a dual polar space $\Delta$ is said to be polarized precisely when all singular hyperplanes of $\Delta$ arise from it, a singular hyperplane of $\Delta$ being formed by the points of $\Delta$ at non maximal distance from a given point, just as in $\AnF$. Thas and Van Maldeghem adopt a different terminology in \cite{TVM}, where the embeddings we call polarized are called full weak embeddings. 
\end{note}

\begin{theo}\label{quot2}
If $\sigma\in \Aut(\FF)\setminus\{\idF\}$ or $\sigma = \idF$ but either $\mathrm{char}(\FF) = 0$ or $\mathrm{char}(\FF)$ is positive and prime to $n+1$, then no proper quotient of $\ve_\sigma$ is polarized. When $\mathrm{char}(\FF)$ is positive and divides $n+1$ then $\veS/[I]$ is the unique polarized proper quotient of $\veS$.  
\end{theo}
{\bf Proof.} Let $S$ be a subspace of $W$ such that $[S]$ defines a polarized proper quotient of $\ve_\sigma$ and let $M$ belong to $S$. By Proposition \ref{prop2}, for every pure tensor $\ba\otimes\alpha\in \Pu((V\otimes V^*)_0)$ the matrix $M$ belongs to the hyperplane $H_{\ba\otimes\alpha^\sigma,\sigma}$ of $W$ spanned by the pure tensors $\bx\otimes\xi\in W$ such that $\xi(\ba)\alpha^\sigma(\bx) = 0$. By Proposition \ref{prop2} we have $H_{\ba\otimes\alpha^\sigma, \sigma} = ((\ba\otimes\alpha^\sigma)^{\perp_f})\cap W$ ($= (\ba\otimes\alpha^\sigma)^{\perp_f}$ when $\sigma \neq \idF$). Hence $M \in H_{\ba\otimes\alpha^\sigma,\sigma}$ if and only if $\ba\otimes\alpha^\sigma \perp_f M$ if and only if $\alpha^\sigma M \ba = 0$. When $\sigma\neq\idF$ the latter condition holds for every $\ba\otimes\alpha\in \Pu((V\otimes V^*)_0)$ if and only if $M = O$. 

Let $\sigma = \idF$. Then $\alpha M \ba = 0$ for every $\ba\otimes\alpha\in \Pu((V\otimes V^*)_0)$ if and only if $M \in \langle I\rangle$. However $M\in W$ and in this case $W = M_{n+1}^0(\FF) = (V\otimes V^*)_0$. Hence $M\neq O$ only if $\mathrm{trace}(I) = 0$, namely $\mathrm{char}(\FF)$ is positive and divides $n+1$.  \hfill $\Box$    

\begin{note}
\em 
The following is implicit in the proof of Theorem \ref{quot2}:  given two non-null matrices $M$ and $N$, the hyperplane $\cH_{N,\sigma}$ arises from $\ve_\sigma/[M]$ if and only if $N\perp_f M$. In particular, assuming that $\veS/[I]$ exists, the hyperplane $\cH_{\ba\otimes\alpha, \idF}$ arises from $\veS/[I]$ if and only $\alpha I\ba = 0$, namely $\alpha(\ba) = 0$. Consequently, the quasi-singular hyperplanes of $\AnF$ which arise from $\veS/[I]$ are precisely the singular ones. 
\end{note} 

\begin{conj}\label{last conj}
Up to isomorphisms, the natural embedding $\veS$, its twistings (when $|\Aut(\FF)| > 1$) and the quotient $\veS/[I]$ (when $\mathrm{char}(\FF) > 0$ divides $n+1$) are the unique polarized embeddings of $\AnF$. 
\end{conj}  

As proved by Thas and Van Maldeghem \cite{TVM}, Conjecture \ref{last conj} holds true when $\FF$ is finite and $n = 2$. 

The property of being polarized is preserved when taking covers. Therefore if Conjecture \ref{last conj} holds true then Conjecture \ref{Main Conj} also holds true.


\begin{thebibliography}{99}
\bibitem{BP2} R.J. Blok and A. Pasini. Point-line geometries with a generating set that depends on the underlying field, in {\em Finite Geometries}  (eds. A. Blokhuis et al.), Kluwer, Dordrecth (2001), 1--25.
\bibitem{BP1} R.J.  Blok and A. Pasini. On absolutely universal embeddings, {\em Dicrete Math.} {\bf 267} (2003), 45-62. 
\bibitem{DSSVM} A. De Schepper, J. Schillewaert and H. Van Maldeghem. On the generatig rank and embedding rank of the hexagonic Lie incidence geometries, submitted. 
\bibitem{FF} C.A. Faure and A. Fr\"{o}licher. {\em Modern Projective Geometry}, Kluwer, Dordrecht 2000. 
\bibitem{Pas} A. Pasini. Geometric hyperplanes of the Lie geometry $\AnF$, Arxiv. preprint 230603947.   
\bibitem{Ron} M.A. Ronan. Embeddings and hyperplanes of discrete geometries, {\em European J. Combin.} {\bf 8} (1987), 179-185.  
\bibitem{S93} E.E. Shult. Embeddings and hyperplanes of Lie incidence geometries, in {\em Groups of Lie Type and Their Geometries} (W.M. Kantor and L. Di martino eds.), {\em L.M.S. Lecture Notes Series} {\bf 207}, Cmabridge Univ. Press, Cambridge (1995), 215-232.  
\bibitem{S} E.E. Shult. {\em Points and Lines}, Springer, Berlin 2011. 
\bibitem{TVM} J. A. Thas and H. Van Maldeghem. Classification of embeddings of the flag geometries of projective planes in finite projective planes, {\em J. Combin. Th. A} {\bf 90} (2000), 159-172 (Part I), 241-256 (Part II), 173-196 (Part III).
\bibitem{V} H. V\"{o}lklein. On the geometry of the adoint representation of a Chevalley group, {\em J. Algebra} {\bf 127} (1989), 139-154.  
\end{thebibliography}
\end{document}